\documentclass[11pt,reqno]{article}

\usepackage{amsfonts}
\usepackage{amscd}
\usepackage{color}
\usepackage{amsmath}
\usepackage{float}
\usepackage{amsthm}
\usepackage{amssymb}
\usepackage{mathrsfs}
\usepackage{amstext}
\usepackage{graphicx}
\usepackage{tikz}
\usepackage{verbatim}
\numberwithin{equation}{section}

\theoremstyle{plain}
\newtheorem{satz}{Theorem}[section]
\newtheorem{defi}[satz]{Definition}
\newtheorem{cor}[satz]{Corollary}
\newtheorem{lem}[satz]{Lemma}
\newtheorem{prop}[satz]{Proposition}
\newtheorem{rem}[satz]{Remark}

\newcommand{\mix}{{\rm mix}}
\newcommand{\re}{\ensuremath{\mathbb{R}}}\newcommand{\N}{\ensuremath{\mathbb{N}}
}
\newcommand{\zz}{\ensuremath{\mathbb{Z}}}\newcommand{\C}{\ensuremath{\mathbb{C}}
}
\newcommand{\T}{\ensuremath{\mathbb{T}^d}}\newcommand{\tor}{\ensuremath{\mathbb{
T}}}

\newcommand{\Z}{{\ensuremath{\zz}^d}}

\newcommand{\n}{\ensuremath{{\N}_0}}
\newcommand{\R}{\ensuremath{{\re}^d}}
\newcommand{\D}{\ensuremath{\mathbb{D}}}
\newcommand\dint{{\rm d}}

\newcommand{\supp}{{\rm supp \, }}

\newcommand{\eps}{\varepsilon}
\newcommand{\rd}{\mathrm{d}}

\newcommand{\bproof}{\begin{proof}}
\newcommand{\eproof}{\end{proof}}

\newlength{\fixboxwidth}
\newcommand{\be}{\begin{equation}}
\newcommand{\ee}{\end{equation}}
\newcommand{\beq}{\begin{eqnarray}}
\newcommand{\beqq}{\begin{eqnarray*}}
\newcommand{\eeq}{\end{eqnarray}}
\newcommand{\eeqq}{\end{eqnarray*}}

\def\Int{\mbox{QMC}}
\def\P{\mathcal{P}}

\def\F{\mathbb{F}}

\begin{document}
\title{Optimal quasi-Monte Carlo rules on order $2$ digital nets for the
numerical integration of multivariate periodic functions}

\author{Aicke Hinrichs$^a$, Lev Markhasin$^b$, Jens Oettershagen$^c$, Tino
Ullrich$^c$\footnote{Corresponding author, Email:
tino.ullrich@hcm.uni-bonn.de}\\\\
$^a$Institute for Analysis, Johannes Kepler University, A-4040 Linz, Austria\\
$^b$Institute for Stochastics and Applications, University of Stuttgart, 70569
Stuttgart\\
$^c$Institute for Numerical Simulation, University of Bonn, 53115 Bonn}

\date{\today}

\maketitle

\begin{abstract}
We investigate quasi-Monte Carlo rules for the numerical integration of
multivariate periodic functions 
from Besov spaces $S^r_{p,q}B(\tor^d)$ with dominating mixed smoothness
$1/p<r<2$. 
We show that order 2 digital nets achieve the optimal rate of convergence
$N^{-r} (\log N)^{(d-1)(1-1/q)}$.
The logarithmic term does not depend on $r$ and hence improves the known bound
of Dick \cite{Di07} for the special case of Sobolev spaces
$H^r_{\text{mix}}(\tor^d)$. Secondly, the rate of convergence is independent of
the integrability $p$ of the Besov space, which allows for sacrificing
integrability in order to gain Besov regularity. Our method combines
characterizations of periodic Besov spaces with dominating mixed smoothness via
Faber bases with sharp estimates of Haar coefficients for the discrepancy
function of order $2$ digital nets. Moreover, we provide numerical computations
which indicate that this bound also holds for the case $r=2$.
\end{abstract}



\section{Introduction}


Quasi-Monte Carlo methods play an important
role for the efficient numerical integration of multivariate functions. Many
real world problems, 
for instance, from finance, quantum physics, meteorology, etc.,
require the computation of integrals of $d$-variate functions where $d$ may be
very large. This can almost never be done analytically since often the available
information of the signal or function $f$ is highly incomplete or simply no
closed-form solution exists. A quasi-Monte Carlo rule approximates the integral
$I(f) = \int_{[0,1]^d} f(x) \, \rd x$
by (deterministically) averaging over $N$ function values taken at fixed points
$X_N = \{x^1,\ldots,x^N\}$, i.e., 
$$
  I_N(X_N,f) := \frac{1}{N}\sum\limits_{i=1}^N f(x^i),
$$
where the $d$-variate function $f$ is assumed to belong to some (quasi-)normed
function space $F_d \subset C([0,1]^d)$. Since the integration weights
$\frac{1}{N}$ are positive and sum up to $1$, QMC integration is stable and easy
to implement which significantly contributed to its popularity.
The QMC-optimal worst-case error with respect to the class $F_d$ is given by 
\be\label{f000}
  \Int_N(F_d):=\inf\limits_{\substack{X_N \subset [0,1]^d\\ \# X_N =
N}}\sup\limits_{\|f|F_d\| \leq 1}
|I(f)-I_N(X_N,f)|\,.
\ee
In this paper we investigate the asymptotical properties of Dick's construction
\cite{Di07} of order $\sigma$ digital nets $X_N$ where $N = 2^n$. This
construction has recently attracted much attention in the area of uncertainty
quantification \cite{SchwabA, SchwabB}. In the present paper we are interested
in the asymptotic optimality of those higher order nets in the sense of
\eqref{f000} with respect to $F_d$ being a periodic Nikol'skij-Besov space
$S^r_{p,q}B(\tor^d)$ with smoothness $r$ larger than $1$ and less than $2$. 

Dick \cite{Di07} showed for periodic Sobolev spaces $H^r_{\text{mix}}(\tor^d) =
S^r_{2,2}B(\tor^d)$\footnote{These spaces are sometimes also referred to as
\emph{Korobov spaces}.}
\begin{equation}\label{f04_2}
  \Int_N(H^r_{\text{mix}}(\tor^d)) \lesssim N^{-\lfloor r \rfloor} (\log
N)^{d\lfloor r \rfloor -1}\quad,\quad N\geq 2\,,
\end{equation}
if $1/2 < r\leq \sigma$. He also considered non-periodic integrands, see
\cite{Di08}. However, well-known asymptotically optimal results for the
integration of periodic Sobolev functions, see for instance the survey
\cite{Te03}, show that the exponent of the $\log$ should be independent of the
smoothness parameter $r$, namely $(d-1)/2$. In that sense, \eqref{f04_2} is far
from being optimal. Nevertheless, Dick's bound \eqref{f04_2} beats the
well-known sparse grid bound if $r$ is an integer and $d$ is large. The latter
bound involves the log-term $(\log N)^{(d-1)(r+1/2)}$, see \cite{DiUl14, Tr10}
and \eqref{f04} below, which represents the best possible rate among all
cubature formulas taking function values on a sparse grid \cite{DiUl14}.

The aim of this paper is twofold. On the one hand we aim at showing the sharp
relation
\begin{equation}\label{f04_10}
  \Int_N(H^r_{\text{mix}}(\tor^d)) \asymp N^{-r} (\log N)^{(d-1)/2}\quad,\quad
N\geq 2\,,
\end{equation}
if $1/2<r<2$ by proving the asymptotical optimality of order $2$ digital nets
for \eqref{f000}. On the other hand we would like to extend \eqref{f04_10} to
periodic Nikol'skij-Besov spaces with dominating mixed smoothness
$S^r_{p,q}B(\tor^d)$, namely, 
\begin{equation}\label{f04_3}
  \Int_N(S^r_{p,q}B(\tor^d)) \asymp N^{- r} (\log N)^{(d-1)(1-1/q)}\quad,\quad
N\geq 2\,,
\end{equation}
for $1/p<r<2$, see Definition \ref{d1} below. An immediate feature of these
error bounds is the fact that the $\log$-term disappears in case $q=1$. Besov
regularity is the correct framework when it comes to integrands of the
form
\be\label{f0002}
	f(x) = \max\{0,g(x)\}\quad,\quad x\in \R\,,
\ee
so-called kink functions, which often occur in mathematical finance, e.g.\ the
pricing of a European call option, whose pay-off function possesses a kink at
the strike price, see e.g. \cite[Chapter 1]{glasserman2004monte}. In general,
one can not expect Sobolev regularity higher than $r=3/2$. However, when
considering Besov regularity we can achieve smoothness $r=2$. Indeed, the simple
example $f(t) = \max\{0,t-1/2\}$ belongs to $B^{2}_{1,\infty}([0,1])$ while its
Sobolev regularity $H^s$ is below $s=3/2$. In a sense, one sacrifices
integrability for gaining regularity. Looking at the bound \eqref{f04_3} above,
we see that cubature methods based on order $2$ digital nets benefit from higher
Besov regularity while the integrability $p$ does not enter the picture. 

Apart from that, spaces of this type have a long history in the former Soviet
Union, see \cite{Am76, Nik75, ScTr87, Te93} and the references
therein. The scale of spaces $S^r_{p,q}B(\tor^d)$ contains several important
special cases of spaces with mixed smoothness like H\"older-Zygmund spaces
$(p=q=\infty)$, the above mentioned Sobolev spaces $(p=q=2)$ and the classical
Nikol'skij spaces $(q=\infty)$. Note that Sobolev spaces $S^r_pH(\tor^d)$ with
integrability $1<p<\infty$ and $r>0$ are not contained in the Besov scale. They
represent special cases of Triebel-Lizorkin spaces $S^r_{p,q}F(\tor^d)$ if
$q=2$. However, classical embedding theorems allow to reduce the question for
$\Int_N(S^r_pH(\tor^d))$ to  \eqref{f04_3} in the case of ``large'' smoothness
$r> \max\{1/p,1/2\}$, see Corollary \ref{mainH} below. For a complete study of
asymptotical error bounds (including the case of ``small'' smoothness) of
numerical integration in spaces $S^r_pH(\tor^d)$ and $S^r_{p,q}F(\tor^d)$ we
refer to the recent preprint \cite{UlUl14}. See also Remark \ref{remH} below. 

The by now classical research topic of numerically integrating 
periodic functions goes back to the work of Korobov \cite{Ko59}, Hlawka
\cite{Hl62}, and Bakhvalov \cite{Ba63} in the 1960s and was continued later by
Temlyakov, see \cite{Te86,Te90,Te91}, Dubinin \cite{Du, Du2} and 
Skriganov~\cite{Sk94}. See also the survey
articles Temlyakov \cite{Te03} and Novak \cite{No15}. In
particular, Temlyakov \cite{Te86, Te90} used the classical Korobov lattice rules
in order to obtain for $1/2<r\leq 1$
\begin{equation}\label{f04_4}
  N^{-r} (\log N)^{(d-1)} \lesssim \Int_N(S^r_{2,\infty}B(\tor^d)) \lesssim N^{-
r} (\log N)^{(d-1)(r+1/2)}
\end{equation}
as well as 
\begin{equation}\label{f04_5}
  N^{-r} (\log N)^{(d-1)/2} \lesssim \Int_N(H^r_{\text{mix}}(\tor^d)) \lesssim
N^{- r} (\log N)^{(d-1)r}
\end{equation}
for $N\geq 2$. In contrast to the quadrature of univariate functions, where
equidistant point
grids lead to optimal formulas, the multivariate problem is much more involved.
In fact, the choice of proper sets $X_N \subset \tor^d$ of integration nodes in
the $d$-dimensional unit cube is the essence of ``discrepancy theory''
\cite{DP, DiPi10} and connected with deep problems
in number theory, already for $d=2$.

Recently, Triebel \cite{Tr10, Tr12} and, independently, D\~ung \cite{Di12}
brought up the framework of tensor Faber bases for functions of the above type.
The main feature is the fact that the basis coefficients are linear combinations
of function values. The corresponding series expansion is thus extremely useful
for sampling and integration issues. Triebel was actually the first who 
investigated cubature formulas for spaces $S^r_{p,q}B(Q_d)$ of functions on the
unit cube $Q_d:=[0,1]^d$.
By using more general cubature formulas of type \eqref{f041} below (with
non-equal weights) and nodes from
a sparse grid Triebel obtained the two-sided estimate 
\begin{equation}\label{f04}
    N^{-r}(\log N)^{(d-1)(1-1/q)}\lesssim \mbox{Int}_N(S^r_{p,q}B(Q_d)) \lesssim
N^{-r}(\log N)^{(d-1)(r+1-1/q)}
\end{equation}
if $1\leq p,q \leq \infty$ and $1/p<r<1+1/p$. Here, $\mbox{Int}_N$ denotes the
optimal worst-case integration error where one admits general (not only QMC)
cubature formulas of type 
\be\label{f041}
    \Lambda_N(X_N,f):= \sum\limits_{i=1}^N \lambda_i f(x^i)\,.
\ee
In contrast to $S^r_{p,q}B(\tor^d)$, the space $S^r_{p,q}B(Q_d)$ consists of
non-periodic functions on $Q_d:=[0,1]^d$. The questions remain how to close the
gaps in the power of the logarithms in \eqref{f04_4}, \eqref{f04_5}, and
\eqref{f04} and what (if existing) are optimal QMC algorithms? 

This question has partly been answered by the first and second named authors for
a subclass of
$S^r_{p,q}B(Q_d)$ with $1/p<r\leq 1$, namely those functions
$S^r_{p,q}B(Q_d)^{\urcorner}$ with vanishing boundary values on the ``upper and
right'' boundary faces, by showing that the lower
bound in \eqref{f04} is sharp for quasi-Monte Carlo methods based on
Chen-Skriganov points, see \cite{Hi10, Ma13_2, Ma13_1,MaDiss13}. Furthermore,
together with M.\ Ullrich the last named author recently observed, that the
classical Frolov method is optimal in all (reasonable) spaces
$S^{r}_{p,q}B(Q_d)^{\square}$ and 
$S^{r}_{p,q}F(Q_d)^{\square}$ of functions with homogeneous boundary condition,
see \cite{UlUl14}. Note, that Frolov's method is an equal-weight cubature
formula of type \eqref{f041} with nodes from a lattice (not a lattice rule). In
a strict sense, Frolov's method is not a QMC method since the weights
$\lambda_i$ do not sum up to $1$.

In this paper we investigate special QMC methods for periodic Nikol'skij-Besov
spaces on $\tor^d$ and answer the above question partly. The picture is clear in
case $d=2$, i.e., for spaces on the $2$-torus $\tor^2$. In fact, we know that 
the lower bound in \eqref{f04} is even sharp for all $r>1/p$, see \cite{Te91,
Ul12_3, DiUl14}. The optimal QMC rule in case $1/p<r<2$ is 
based on Hammersley points \cite{Ul12_3}, which provide the optimal discrepancy
in this setting \cite{Hi10}. This paper can be seen as continuation of
\cite{Ul12_3} for the higher-dimensional situation by adopting methods from
\cite{Hi10, Ma13_1, Ma13_2, Ma15}. 

We will prove the optimality of QMC methods based on order $2$ digital nets in
the framework of periodic Besov 
spaces with dominating mixed smoothness if $1/p<r<2$. Due to the piecewise
linear building blocks, we can not expect to
get beyond $1/p<r<2$ with our proof method even when taking higher order digital
nets. Therefore, 
this restriction seems to be technical and may be overcome by using smoother 
basis atoms like piecewise polynomial B-splines \cite{Di12}.

We illustrate our theoretical results with numerical computations in the Hilbert
space case $H^r_{\mix}(\tor^d)$
in several dimensions $d$ and for different smoothness parameters $r$. In the
case of integer smoothness we exploit 
an exact representation formula for the worst-case integration error of an
arbitrary cubature rule. A numerical evaluation of this formula 
indicates that the results in Theorem \ref{main_theorem} below keep valid for
$r=2$. The comparison with other widely used cubature rules such as sparse grids
and
Halton points in all dimensions and Fibonacci lattices in dimension $d=2$ shows
that order 2 digital nets perform very well
not only asymptotically but already for a relatively small number of sample
points. Finally, we consider a simple test function which is a tensor product of
univariate functions of the form \eqref{f0002}. Expressing the regularity of
such functions in Besov spaces of mixed smoothness allows the correct prediction
of the asymptotical rate of the numerical integration error which is verified by
our numerical experiments. However, the applicability of order $2$ nets to
real-world problems from option pricing etc., where the kinks are not
necessarily axis aligned, requires further research.

The paper is organized as follows. In Section 2 we introduce the function spaces
of interest and provide the necessary characterizations and properties. The
classical definition by mixed iterated differences will turn out to be of
crucial importance. In Section 3 we deal with the Faber and Haar basis,
especially with their hyperbolic (anisotropic) tensor product. The main tools
represent Propositions \ref{discvscont} and \ref{contvsdisc} where the function
space norm is related to the Faber coefficient sequence space norm and vice
versa. In Section $4$ we recall Dick's construction of higher order digital nets
and compute the Haar coefficients of the associated discrepancy function. We
continue in Section $5$ by interpreting the Haar coefficients of the discrepancy
function in terms of integration errors for tensorized Faber hat functions.
Combining those estimates with the Faber basis expansion and the
characterization from Section $3$ we obtain our main results in Theorem
\ref{main_theorem}. Finally, 
Section $6$ provides the numerical results.

{\bf Notation.} As usual $\N$ denotes the natural numbers, $\N_0=\N\cup\{0\}$,
$\N_{-1}=\N_0\cup\{-1\}$, $\zz$ denotes the integers, 
$\re$ the real numbers, and $\C$ the complex numbers. The letter $d$ is always
reserved for the underlying dimension in $\re^d, \zz^d$ etc. We denote
by $\langle x,y\rangle$ the usual Euclidean inner product and inner products in
general. For $a\in \re$ we denote $a_+ := \max\{a,0\}$. 
For $0<p\leq \infty$ and $x\in \R$ we denote $|x|_p = (\sum_{i=1}^d
|x_i|^p)^{1/p}$ with the
usual modification in the case $p=\infty$. We further denote 
$x_+ := ((x_1)_+,\ldots,(x_d)_+)$ and $|x|_+ := |x_+|_1$. By $(x_1,\ldots,x_d)>0$
we mean that each coordinate is positive. By $\tor$ we denote the torus
represented by the interval $[0,1]$, where the end points are identified.
If $X$ and $Y$ are two (quasi-)normed spaces, the (quasi-)norm
of an element $x$ in $X$ will be denoted by $\|x|X\|$. 
The symbol $X \hookrightarrow Y$ indicates that the
identity operator is continuous. For two sequences $a_n$ and $b_n$ we will write
$a_n \lesssim b_n$ if there exists a constant $c>0$ such that $a_n \leq c\,b_n$
for all $n$. We will write $a_n \asymp b_n$ if $a_n \lesssim b_n$ and $b_n
\lesssim a_n$.


\section{Periodic Besov spaces with dominating mixed smoothness}
\label{perfs}
Let $\tor^d$ denote the $d$-torus, represented in the Euclidean space $\re^d$
by the cube $\tor^d = [0,1]^d$, where opposite points are identified.
That means $x,y\in \re^d$ are identified if and only if $x-y = k$, where $k =
(k_1,\ldots,k_d) \in \zz^d$. 
The computation of the Fourier coefficients of an integrable $d$-variate
periodic function is performed by the formula 
$$
   \hat{f}(k)= \int_{\tor^d}
  f(x)e^{-i2\pi k \cdot x}\,\dint x\quad,\quad k\in \zz^d\,.
$$
Let further denote $L_p(\tor^d)$, $0<p\leq \infty$, the space of
all measurable functions $f:\tor^d\rightarrow \C$ satisfying
$$
    \|f\|_p =
    \Big(\,\int_{\tor^d}|f(x)|^p\,\dint x\Big)^{1/p} < \infty
$$
with the usual modification in case $p=\infty$. The space $C(\tor^d)$ is often
used as a replacement for $L_{\infty}(\tor^d)$. It denotes the collection of all
continuous and
bounded periodic functions equipped with the $L_{\infty}$-topology. 

\subsection{Definition and basic properties}

In this section we give the definition of Besov spaces with dominating mixed
smoothness on $\tor^d$ based on a dyadic decomposition on the Fourier side. We
closely follow \cite[Chapter 2]{ScTr87}.
To begin with, we recall the concept of a dyadic decomposition of unity. The
space $C^{\infty}_0(\re^d)$ consists of all infinitely many times differentiable
compactly supported functions. 
    \begin{defi}\label{cunity} Let $\Phi(\re)$ be the collection of all
      systems
       $\varphi = \{\varphi_n(x)\}_{n=0}^{\infty} \subset C^{\infty}_0(\re^d)$
       satisfying
       \begin{description}
       \item(i) ${\supp}\,\varphi_0 \subset \{x:|x| \leq 2\}$\, ,
       \item(ii) ${\supp}\,\varphi_n \subset \{x:2^{n-1} \leq |x|
       \leq 2^{n+1}\}\quad,\quad n= 1,2,\ldots ,$
       \item(iii) For all $\ell \in \N_0$ it holds
       $\sup\limits_{x,n}
       2^{n\ell}\, |D^{\ell}\varphi_n(x)| \leq c_{\ell} <\infty$\, ,
       \item(iv) $\sum\limits_{n=0}^{\infty} \varphi_n(x) = 1$ for all
       $x\in \re$.
       \end{description}
    \end{defi}

    \begin{rem}\label{speciald}
      The class $\Phi(\re)$ is not empty. We consider the following
      standard example.
      Let $\varphi_0(x)\in S(\re)$ be a smooth function with $\varphi_0(x) =
      1$ on $[-1,1]$ and $\varphi_0(x) = 0$
      if $|x|>2$. For $n>0$ we define
      $$
         \varphi_n(x) = \varphi_0(2^{-n}x)-\varphi_0(2^{-n+1}x).
      $$      It is easy to verify that the system $\varphi =
      \{\varphi_n(x)\}_{n=0}^{\infty}$ satisfies (i) - (iv).
    \end{rem}
    \noindent Now we fix a system $\varphi=\{\varphi_n\}_{n\in \zz} \in
\Phi(\re)$, where we 
    put $\varphi_n \equiv 0$ if $n<0$. For $j = (j_1,\ldots,j_d) \in \zz^d$ let the
building blocks $f_{j}$ be given by
    \begin{equation}\label{f2}
	f_j(x) = \sum\limits_{k\in \zz^d}
	\varphi_{j_1}(k_1)\cdots\varphi_{j_d}(k_d)\hat{f}(k)e^{i 2\pi
k\cdot x}\quad,\quad x\in \tor^d\,,j\in \zz^d\,.
    \end{equation}

    \begin{defi}\label{d1} (Mixed periodic Besov and Sobolev space) \\
    {\em (i)}
    Let $0< p,q\leq
    \infty$ and $r>\sigma_p:=(1/p-1)_+$. Then $S^{r}_{p,q}B(\tor^d)$ is defined
as the collection of all $f\in
    L_1(\tor^d)$ such that
    \begin{equation}\label{f3}
	\|f|S^r_{p,q}B(\tor^d)\|^{\varphi}:=\Big(\sum\limits_{j\in \N_0^d}
2^{|j|_1rq}\|f_j\|_p^q\Big)^{1/q}
    \end{equation}
    is finite (usual modification in case $q=\infty$).\\
    {\em (ii)} Let $1<p<\infty$ and $r>0$. Then $S^r_pH(\tor^d)$ is defined as
the collection of all $f\in L_p(\tor^d)$ such that 
    $$
	\|f|S^r_pH(\tor^d)\|^{\varphi}:=\Big\|\Big(\sum\limits_{j\in
\N_0^d}2^{|j_1|r2}|f_j(x)|^2\Big)^{1/2}\Big\|_p
    $$
    is finite. 
    \end{defi}
Recall, that this definition is independent of the chosen system $\varphi$ in
the sense of equivalent (quasi-)norms. Moreover, in case $\min\{p,q\}\geq 1$ the
defined spaces are Banach spaces, whereas they are quasi-Banach spaces in case
$\min\{p,q\} < 1$. For details confer \cite[2.2.4]{ScTr87}.
In this paper we are mainly concerned with spaces providing sufficiently large
smoothness $(r>1/p)$ such that the elements (equivalence classes) in
$S^r_{p,q}B(\tor^d)$ contain a continuous representative. We have the
following elementary embeddings, see \cite[2.2.3]{ScTr87}.

\begin{lem}\label{emb} Let $0<p<\infty$, $r\in \re$, and $0<q\leq \infty$.
\begin{description}
   \item (i) If $\varepsilon>0$ and $0<v\leq \infty$ then
 $$
   S^{r+\varepsilon}_{p,q}B(\tor^d) \hookrightarrow S^r_{p,v}B(\tor^d)\,.
 $$
 \item (ii) If $p<u\leq \infty$ and $r-1/p = t-1/u$ then
 $$
     S^r_{p,q}B(\tor^d) \hookrightarrow S^t_{u,q}B(\tor^d)\,.
 $$
 \item (iii) If $r>1/p$ then
 $$
    S^r_{p,q}B(\tor^d) \hookrightarrow C(\tor^d)\,.
 $$
 \item(iv) If $1<p<\infty$ and $r>0$ then
 $$
    S^r_{p,\min\{p,2\}}B(\tor^d) \hookrightarrow S^r_pH(\tor^d) \hookrightarrow
S^r_{p,\max\{p,2\}}B\,.
 $$
 \item (v) If $2\leq p<\infty$ and $r>0$ then 
 $$
    S^r_pH(\tor^d) \hookrightarrow S^r_2H(\tor^d) = H^r_{\text{mix}}(\tor^d) =
S^r_{2,2}B(\tor^d)\,.
 $$
\end{description}
\end{lem}

\subsection{Characterization by mixed differences}

In this subsection we will provide the classical characterization by mixed
iterated differences as it is used for instance in \cite{Am76}. The main issue
will be the equivalence of both approaches, the Fourier analytical approach in
Definition \ref{d1} and the difference approach, see Lemma \ref{diff} below. We
will need some tools from Harmonic Analysis, the Peetre maximal function and the
associated maximal inequality, see \cite[1.6.4, 3.3.5]{ScTr87}. 
For $a>0$ and $b = (b_1,\ldots,b_d)>0$ we define
the Peetre maximal function $P_{b,a}f$ for a trigonometric polynomial $f$,
i.e., 
$$
      P_{b,a}f(x) := \sup\limits_{y\in \R} \frac{|f(y)|}{(1+b_1|x_1-y_1|)^a\cdots (1+b_d|x_d-y_d|)^a}
$$
The following maximal inequality for multivariate trigonometric polynomials with
frequencies in the rectangle
$Q_b:=[-b_1,b_1]\times \ldots \times [-b_d,b_d]$ will be of crucial importance.
\begin{lem}\label{peetremaximal} Let $0<p\leq \infty$, $b = (b_1,\ldots,b_d)>0$,
and $a>1/p$. Let further 
$$
  f = \sum_{\substack{|k_i| \leq b_i\\ i = 1,\ldots,d}} \hat{f}(k)e^{2\pi ik\cdot
x}
$$
be a trigonometric polynomial with frequencies in the rectangle $Q_b$\,. 
Then a constant $c>0$ independent of $f$ and $b$ exists such that 
$$
    \|P_{b,a}f\|_p \leq c\|f\|_p\,. 
$$
\end{lem}
\noindent Now we introduce the basic concepts of iterated differences
$\Delta^m_h(f,x)$ of a function $f$. For univariate functions $f:\tor \to \C$
the $m$th difference operator $\Delta_h^{m}$ is defined by
\begin{equation*}
\Delta_h^{m}(f,x) := \sum_{j =0}^{m} (-1)^{m - j} \binom{m}{j} f(x +
jh)\quad,\quad x\in \tor, h\in [0,1]\,.
\end{equation*}

The following Lemma states an important relation between $m$th order differences
and Peetre maximal functions of trigonometric polynomials, see \cite[Lemma
3.3.1]{Ul06}.

\begin{lem}\label{ddim} Let $a,b>0$ and 
$$
    f = \sum\limits_{|k|\leq b} \hat{f}(k)e^{2\pi ikx}\quad,\quad x\in \tor\,,
$$
be a univariate trigonometric polynomial with frequencies in $[-b,b]$. Then
there exists a constant $c>0$ independent of $f$ such that for every $h\in \re$
  \begin{equation}
         |\Delta^{m}_h(f,x)| \leq c\min\{1,|bh|^{m}\}\max\{1,|b
h|^a\}P_{b,a}f(x)\quad,\quad x\in \tor\,.
  \end{equation}
 \end{lem}

In order to characterize multivariate functions we need the concept of mixed
differences with respect to coordinate directions. 
Let $e$ be any subset of
$\{1,\ldots,d\}$. For multivariate functions $f:\tor^d\to \C$ and $h\in [0,1]^d$
the mixed $(m,e)$th difference operator $\Delta_h^{m,e}$ is defined by
\begin{equation*}
\Delta_h^{m,e} := \
\prod_{i \in e} \Delta_{h_i,i}^m\quad\mbox{and}\quad \Delta_h^{m,\emptyset} = 
\operatorname{Id},
\end{equation*}
where $\operatorname{Id}f = f$ and $\Delta_{h_i,i}^m$ is the univariate operator
applied to the $i$-th coordinate of $f$ with the other variables kept fixed. Let
us further define the mixed $(m,e)$th modulus of continuity by 
\begin{equation}\label{modc}
\omega_{m}^e(f,t)_p:= \sup_{|h_i| < t_i, i \in
e}\|\Delta_h^{m,e}(f,\cdot)\|_{p}\quad,\quad t \in [0,1]^d,
\end{equation}
for $f \in L_p(\tor^d)$ (in particular, $\omega_{m}^{\emptyset}(f,t)_p =
\|f\|_{p}$)\,. 
We aim at an equivalent characterization of the Besov spaces
$S^r_{p,q}B(\tor^d)$. The following lemma answers this question partially. There
are 
still some open questions around this topic, see for instance \cite[2.3.4,
Remark 2]{ScTr87}. The following Lemma is a straight-forward 
modification of \cite[Theorem 4.6.1]{Ul06}.

\begin{lem}\label{diff} Let $1\leq p\leq \infty$, $0< q \leq \infty$
and $m\in \N$ with $m>r>0$. Then 

$$
    \|f|S^r_{p,q}B(\tor^d)\|^{\varphi} \asymp
\|f|S^r_{p,q}B(\tor^d)\|^{(m)}\quad,\quad f\in L_1(\tor^d)\,,
$$
where 
\begin{equation}\label{f4}
     \|f|S^r_{p,q}B(\tor^d)\|^{(m)} := \Big[\sum\limits_{j\in \N_{0}^d}
2^{r|j|_1q}\omega_{m}^{e(j)}(f,2^{-j})^q_p\Big]^{1/q}\,.
\end{equation}
In case $q=\infty$ the sum in \eqref{f4} is replaced by the $\sup$ over $j$.
Here, $e(j) = \{i:j_i \neq 0\}$. 

\end{lem}
\bproof This assertion is a modified version of
\cite[Theorem 4.6.2]{Ul06} for the bivariate setting. Let us recall
some basic steps in the proof. The relation
$$
     \|f|S^r_{p,q}B(\tor^d)\|^{(m)} \leq C_1\|f|S^r_{p,q}B(\tor^d)\|^{\varphi}
$$
is obtained by applying \cite[Lemma 3.3.2]{Ul06} to the building blocks $f_j$
in \eqref{f2},
which are indeed trigonometric polynomials, and using the proof technique in
\cite[Theorem 3.8.1]{Ul06}.

To obtain the converse relation
$$
    \|f|S^r_{p,q}B(\tor^d)\|^{\varphi} \leq C_2\|f|S^r_{p,q}B(\tor^d)\|^{(m)}
$$
we take into account the characterization via rectangle means given in
\cite[Theorem 4.5.1]{Ul06}. We apply the techniques in Proposition
3.6.1 to switch from rectangle means to moduli of smoothness by
following the arguments in the proof of Theorem 3.8.2. It remains to
discretize the outer integral (with respect to the step length of the
differences) in order to replace it by a sum. This is done by standard
arguments. Thus, we almost arrived at \eqref{f4}. Indeed, the final step is to
get rid of those summands where the summation index is negative. But this is
trivially done by omitting the corresponding difference (translation invariance
of $L_p$) such that the respective sum is just a converging geometric series
(recall that $r>0$). \eproof

\begin{rem}\label{intmeans} By replacing the moduli of continuity \eqref{modc}
by more regular variants like integral means 
of differences \cite{Ul06} we can extend the characterization in Lemma
\ref{diff} to all $0<p\leq \infty$ and $r>(1/p-1)_+$, see also 
Remark \ref{rectmeans} below.
\end{rem}

\section{Haar and Faber bases}
\label{HaarFaber}
\begin{figure}[H]\label{fig:3}
\centering
\begin{tikzpicture}[scale=3]

\draw[->] (-0.1,0.0) -- (1.1,0.0);
\draw[->] (0.0,-0.1) -- (0.0,1.1); 

\draw (1.0,0.03) -- (1.0,-0.03) node [below] {$1$};
\draw (0.03,1.0) -- (-0.03,1.00) node [left] {$1$};
\draw[->] (0.8,0.9) -- (0.7,0.7);

\node at (0.8,1) {$v_{0,0}$};

\node at (1.1,0.5) {$j=0$};

\draw (0,0) -- (0.5,1);
\draw (0.5,1) -- (1,0);

\draw[->] (1.4,0.0) -- (2.6,0.0); 
\draw[->] (1.5,-0.1) -- (1.5,1.1); 

\draw (2.5,0.03) -- (2.5,-0.03) node [below] {$1$};
\draw (1.53,1.0) -- (1.47,1.00) node [left] {$1$};
\draw[->] (2.5,0.9) -- (2.4,0.7) ;
\draw[->] (2,0.9) -- (1.9,0.7) ;

\node at (2.5,1) {$v_{1,1}$};
\node at (2,1) {$v_{1,0}$};

\node at (2.6,0.5) {$j=1$};

\draw (1.5,0,0) -- (1.75,1.0);
\draw (1.75,1) -- (2,0.0);
\draw (2,0) -- (2.25,1.0);
\draw (2.25,1) -- (2.5,0.0);

\draw[->] (2.9,0.0) -- (4.1,0.0); 
\draw[->] (3.0,-0.1) -- (3.0,1.1);

\draw (4.0,0.03) -- (4.0,-0.03) node [below] {$1$};
\draw (3.03,1.0) -- (2.97,1.00) node [left] {$1$};

\fill[gray!60] plot[domain=3:3.5] (\x,-6+2*\x)%
                      -- plot[domain=4:3.5] (\x,0);%

\node at (4.3,0.5) {$j\in\{0,1\}$};

\draw (3,0) -- (3.5,1);
\draw (3.5,1) -- (4,0);

\draw (3.0,0,0) -- (3.25,1.0);
\draw (3.25,1) -- (3.5,0.0);
\draw (3.5,0) -- (3.75,1.0);
\draw (3.75,1) -- (4,0.0);

\end{tikzpicture}
  \caption{Univariate hierarchical Faber basis  on $\tor$ for levels $j\in
\{0,1\}$ and their union.} \label{fig_Faber1}
\end{figure}
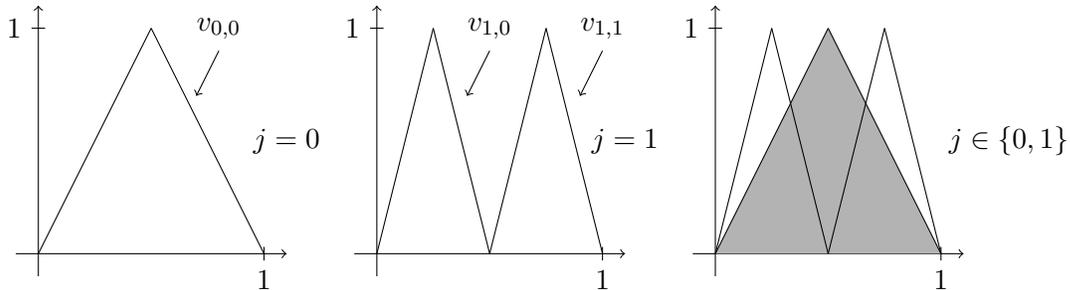

\subsection{The tensor Haar basis}

Let $h(x)$ denote a piecewise constant step function on the real line
given by 
$$
      h(x) = \left\{\begin{array}{rcl}
                    1&:&0\leq x \leq 1/2\,,\\
                    -1&:&1/2<x \leq 1\,,\\
                    0&:& \mbox{otherwise}\,.
                \end{array}\right.
$$
For $j\in \N_0$ and $k \in \D_j := \{0,1,\ldots,2^j-1\}$ we put
$$
      h_{j,k} := h(2^j\cdot-k)\,.
$$
Clearly $h_{j,k}$ is supported in $[0,1]$ for $j\in \N_0$, $k\in \D_j$\,. Let
now $j=(j_1,\ldots,j_d) \in \N_0^d$ and 
$k\in \D_j := \D_{j_1}\times\ldots\times \D_{j_d}$\,. We denote by
$$
    h_{j,k}(x_1,\ldots,x_d) := h_{j_1,k_1}(x_1)\cdots
h_{j_d,k_d}(x_d)\quad,\quad (x_1,\ldots,x_d)\in [0,1]^d\,,
$$ 
the tensor Haar function with respect to the level $j\in \N_0^d$ and the
translation $k\in \Z$
and 
$$
    \mu_{j,k}(f) = \int_{[0,1]^d} f(x)h_{j,k}(x)\,\dint x
$$
the corresponding Haar coefficient for $f\in L_1([0,1]^d)$.

\subsection{The univariate Faber basis}

Recently, Triebel \cite{Tr10,Tr12} and, independently, D\~ung \cite{Di12}
observed
the potential of the Faber basis for the
approximation and integration of functions with dominating mixed smoothness.
The latter reference is even more general and uses so-called B-spline
representations of functions, where the Faber system is a special case. We note
that the Faber basis also plays an important role in the  construction of sparse
grids which go back to \cite{Smolyak:1963} and are utilized in many applications
for the discretization and approximation of function spaces with dominating
mixed smoothness, see e.g. \cite{Bungartz.Griebel:2004, Ul08}.

Let us briefly recall the basic facts about the Faber basis taken from
\cite[3.2.1, 3.2.2]{Tr10}. 
\begin{defi} Let $v(x)$ be the $L_{\infty}$-normalized integrated Haar
function $h$, i.e., 
\be\label{vinth}
  v(x):=2\int_0^x h(t)\,\dint t\quad,\quad x\in \re\,,
\ee
and for $j\in \n, k\in \D_j$  
\be\label{f19}
      v_{j,k} = v(2^j\cdot - k)\,.
\ee
For notational reasons we let $v_{-1,0}:=1$ and obtain the Faber system
$$
    F:=\{v_{j,k}:j\in \N_{-1}, k\in \D_j\}\,.
$$
\end{defi}

Faber \cite{Fa09} observed that every continuous
(non-periodic) function $f$ on $[0,1]$ can be represented as
\begin{equation}\label{f51}
    f(x) = f(0)\cdot (1-x)+f(1)\cdot x -
\frac{1}{2}\sum\limits_{j=0}^{\infty}\sum\limits_{k=0}^{2^j-1}
\Delta^2_{2^{-j-1}}(f,2^{-j}k)v_{j,k}(x)
\end{equation}
with uniform convergence, see e.g. \cite[Theorem 2.1, Step 4]{Tr10}.
Consequently, every periodic function from 
$C(\tor)$ can be decomposed via the ($1$-periodic extended) system $F$ such
that 
\begin{equation}\label{f5}
    f(x) = f(0) - \frac{1}{2}\sum\limits_{j=0}^{\infty}\sum\limits_{k=0}^{2^j-1}
    \Delta^2_{2^{-j-1}}(f,2^{-j}k)v_{j,k}(x)
\end{equation}
with convergence in $C(\tor)$\,.

\subsection{The tensor Faber basis}
Let now $f(x_1,\ldots,x_d)$ be a $d$-variate function $f\in C(\tor^d)$. By fixing
all variables except $x_i$
we obtain by $g(\cdot) = f(x_1,\ldots, x_{i-1},\cdot,x_{i+1},\ldots,x_d)$ a univariate
periodic continuous function. By applying \eqref{f5} in every such component we
obtain the representation in $C(\T)$
\begin{equation}\label{repr}
  f(x) = \sum\limits_{j\in \N_{-1}^d} \sum\limits_{k\in \D_j} d^2_{j,k}(f)
  v_{j,k}(x)\quad,\quad x\in \tor^d\,,
\end{equation}
where 
$$
  v_{j,k}(x_1,\ldots,x_d):=v_{j_1,k_1}(x_1)\cdots
v_{j_d,k_d}(x_d)\quad,\quad j\in \N_{-1}^d, k\in \D_j\,,
$$
and 
\begin{equation}\label{f100}
      d^2_{j,k}(f) =
(-2)^{-|e(j)|}\Delta^{2,e(j)}_{2^{-(j+1)}}(f,x_{j,k})\quad,\quad j\in \N_{-1}^d,
k\in \D_j\,.
\end{equation}
Here we put $e(j) = \{i:j_i \neq -1\}$ and $x_{j,k} =
(2^{-(j_1)_+}k_1,\ldots,2^{-(j_d)_+}k_d)$\,.

Our next goal is to discretize the spaces $S^r_{p,q}B(\tor^d)$ using the Faber
system $\{v_{j,k}\,:\,j\in \N_{-1}^d, k\in \D_j\}$. We obtain a sequence space
isomorphism performed by the coefficient
mapping $d^2_{j,k}(f)$ above. In \cite[3.2.3,
3.2.4]{Tr10} and \cite[Theorem 4.1]{Di12} this was done for the non-periodic
setting $S^r_{p,q}B(Q_d)$. Our
proof is completely different and uses only classical tools.
This makes the proof a bit more transparent and
self-contained. With these tools we show that one direction of the
equivalence relation can be extended to $1/p<r<2$.

\begin{defi}\label{defsequ} Let $0<p,q\leq \infty$ and $r\in \re$. Then
$s^r_{p,q}b$ is
the collection of all sequences $\{\lambda_{j,k}\}_{j\in \N_{-1}^d, k\in \D_j}$
  such that
$$
  \|\lambda_{j,k}|s^r_{p,q}b\|:=\Big[\sum\limits_{j\in
\N_{-1}^d}2^{|j|_1(r-1/p)q}\Big(\sum\limits_{k\in
\D_j}|\lambda_{j,k}|^p\Big)^{q/p} \Big]^{1/q}
$$
is finite.
\end{defi}

\begin{lem}\label{qBanach} Let $0<p,q\leq \infty$ and $r\in \re$. The
space $s^r_{p,q}b$ is a Banach space if $\min\{p,q\}\geq 1$. In case
$\min\{p,q\} <1$ the space $s^r_{p,q}b$ is a quasi-Banach space. Moreover, if
$u:=\min\{p,q,1\}$ it is a $u$-Banach space, i.e.,
$$
    \|\lambda+\mu|s^r_{p,q}b\|^u \leq \|\lambda|s^r_{p,q}b\|^u +
\|\mu|s^r_{p,q}b\|^u\quad,\quad \lambda,\mu \in s^r_{p,q}b\,.
$$
\end{lem}

\begin{prop}\label{discvscont} Let $1/2 < p \leq \infty$, $0<q\leq \infty$ and
$1/p<r<2$.
Then there
exists a constant $c>0$ such that
\begin{equation}\label{f7}
    \big\|d^2_{j,k}(f)|s^r_{p,q}b\big\| \leq c
\|f|S^r_{p,q}B(\tor^d)\|^{\varphi}
\end{equation}
for all $f\in C(\tor^d)$.

\end{prop}
\bproof The main idea is the same as in the proof of Lemma
\ref{diff}. We make use of the decomposition \eqref{f2} in a slightly modified
way. Let us first
assume $1\leq p,q\leq \infty$. The modifications in case $\min\{p,q\}<1$ are
straight-forward. For fixed $j\in
\N_{-1}^d$ we write $f = \sum_{\ell \in \zz^d} f_{j+\ell}$. Putting this
into
\eqref{f7} and using the triangle inequality yields
\begin{equation}\label{f8}
  \begin{split}
   \big\|d^2_{j,k}(f)|s^r_{p,q}b\big\| &\asymp 
   \Big[\sum\limits_{j\in
   \N_{-1}^d}2^{r|j|_1q}\Big\|\sum\limits_{k\in \D_j}
d^2_{j,k}(f)v_{j,k}(x)\Big\|_p^q \Big]^{1/q}\\
   &\leq \sum\limits_{\ell \in \zz^d}
   \Big[\sum\limits_{j\in
   \N_{-1}^d}2^{r|j|_1q}\Big\|\sum\limits_{k\in \D_j}
d^2_{j,k}(f_{j+\ell})v_{j,k}(x)\Big\|_p^q\Big]^{1/q}\,.
  \end{split}
\end{equation}
Let us continue in deriving a point-wise upper bound for the absolute value of
the function 
$$
   F_{j,\ell}(x):=\sum\limits_{k\in \D_j} d^2_{j,k}(f_{j+\ell})v_{j,k}(x)\,.
$$
Clearly, we have
\be\label{peetre0}
    |F_{j,\ell}(x)| \leq |d^2_{j,k}(f_{j+\ell})| \lesssim
|\Delta^{2,e(j)}_{2^{-(j+1)}}(f_{j+\ell},x_{j,k})|
\ee
whenever $x \in [2^{-(j_1)_+}k_1,2^{-(j_1)_+}(k_1+1)] \times\ldots\times
[2^{-(j_d)_+}k_d,2^{-(j_d)_+}(k_d+1)]$\,. Let us estimate the iterated
differences $\Delta^{2,e(j)}_{2^{-(j+1)}}(f_{j+\ell},x_{j,k})$ one by one. For
$i\in e(j)$ we have for $x$ such that 
$|x_{j_i,k_i}-x| \leq 2^{-j_i}$ the bound 
$$
      |\Delta^2_{2^{-(j_i+1)}}(g,x_{j_i,k_i})| \lesssim \sup\limits_{|y|\lesssim
2^{-j_i}} |g(x+y)| \lesssim \sup\limits_{|y|\lesssim 2^{-j_i}}
\frac{|g(x+y)|}{(1+2^{j_i}|y|)^a} \leq P_{2^{j_i},a}g(x)
$$
for a univariate continuous function $g$. In case $i \notin e(j)$ we have $j_i =
-1$ and 
$$
    |g(0)| \leq \sup\limits_{|y|\leq 1} |g(x+y)| \lesssim \sup\limits_{|y|\leq
1} \frac{|g(x+y)|}{(1+2^{j_i}|y|)^a} \leq P_{2^{j_i},a}g(x)\,.
$$
If $\ell_i \geq 0$ then, by definition,  
\be\label{peetre1}
  P_{2^{j_i},a}g(x) \leq 2^{\ell_i a} P_{2^{j_i+\ell_i},a}g(x)\,.
\ee
On the other hand, the estimate in Lemma \ref{ddim} gives in case $i\in e(j)$
for a univariate trigonometric polynomial $g_{j_i+\ell_i}(t) = \sum_{k\in
\zz}\varphi_{j_i+\ell_i}(k)\hat{g}(k)e^{2\pi i kt}$
$$
  |\Delta^2_{2^{-(j_i+1)}}(g_{j_i+\ell_i},x_{j_i,k_i})| \lesssim
\min\{1,2^{2\ell_i}\}\max\{1,2^{\ell_i
a}\}P_{2^{\ell_i+j_i},a}g_{j_i+\ell_i}(x_{j_i,k_i})\,.   
$$
If $\ell_i<0$ and $|x-x_{j_i,k_i}| \leq 2^{-j_i}$ this reduces to 
\be\label{peetre2}
    |\Delta^2_{2^{-(j_i+1)}}(g_{j_i+\ell_i},x_{j_i,k_i})| \lesssim
2^{2\ell_i}P_{2^{\ell_i+j_i},a}g_{j_i+\ell_i}(x)\,.
\ee	
Note, that in case $i\notin e(j)$ there is nothing to prove since
$g_{j_i+\ell_i} \equiv 0$. 
Applying the point-wise estimates \eqref{peetre1} and \eqref{peetre2} to the
right-hand side of \eqref{peetre0} we obtain
$$
    |F_{j,\ell}(x)| \lesssim P_{2^{\ell+j},a}f_{j+\ell}(x)\prod\limits_{i\in
e(j)} \min\{2^{2\ell_i},1\}\max\{2^{\ell_ia},1\}\,,
$$
where $2^{\ell+j}:=(2^{\ell_1+j_1},\ldots,2^{\ell_d+j_d})$\,. Using the Peetre
maximal inequality, Lemma \ref{peetremaximal} yields
\begin{equation}
  \begin{split}
    \|F_{j,\ell}\|_p &\lesssim \|P_{2^{\ell+j},a}f_{j+\ell}\|_p\cdot
\prod\limits_{i\in e(j)} \min\{2^{2\ell_i},1\}\max\{2^{\ell_ia},1\}\\
    &\lesssim \|f_{j+\ell}\|_p\prod\limits_{i\in e(j)}
\min\{2^{2\ell_i},1\}\max\{2^{\ell_ia},1\}\,,
  \end{split}  
  \end{equation}
whenever $a>1/p$\,. If $r>1/p$ we can choose 
\be\label{condi}
    1/p<a<r<2\,.
\ee
Therefore, if $\ell \in \zz^d$ 
\be\label{anorm2}
    \sum\limits_{j\in \N_{-1}^d} 2^{r|j|_1 q} \|F_{j,\ell}\|_p^q \lesssim 
    \sum\limits_{j\in \N_{-1}^d}2^{r|j+\ell|_1 q}\|f_{j+\ell}\|_p^q
\prod\limits_{i=1}^d A^q_{\ell_i}\,,
\ee
where for $n\in \zz$
\be\label{Amn}
    A_{n} = \left\{\begin{array}{rcl}
		    2^{(2-r)n}&:&n<0,\\
		    2^{(a-r)n}&:&n\geq 0.
		  \end{array}\right.
\ee
Under the condition \eqref{condi} it follows from \eqref{Amn} that there is a
$\delta>0$ such that $A_{n} \leq 2^{-\delta |n|}$ and hence
$$
    \sum\limits_{j\in \N_{-1}^d} 2^{r|j|_1 q} \|F_{j,\ell}\|_p^q \lesssim
2^{-q\delta|\ell|_1}\|f|S^r_{p,q}B(\tor^d)\|^q\cdot 
$$
Plugging this into \eqref{f8} yields \eqref{f7}\,.\eproof

Let us prove the converse statement. The version below slightly differs from its
$2$-dimensional counterpart given in \cite{Ul12_3} although the proof technique
is the same. We observed that the restriction $r>1/p$ is actually not required. 

\begin{prop}\label{contvsdisc} Let $1\leq p \leq \infty$, $0<q\leq \infty$,
$0<r<1+1/p$. Then there exists a constant $c>0$ such that
\begin{equation}\label{f52}
    \|f|S^r_{p,q}B(\tor^d)\|^{\varphi} \leq c\|d^2_{j,k}(f)|s^r_{p,q}b\|
\end{equation}
for all $f\in C(\tor^d)$ with finite right-hand side \eqref{f52}\,.

\end{prop}

\bproof We use the characterization in Lemma \ref{diff}, which says that
$$
   \|f|S^r_{p,q}B(\tor^d)\|^{\varphi} \lesssim \|f|S^r_{p,q}B(\tor^d)\|^{(m)}
$$
for some fixed $m\geq 2$\,. Let us assume $1\leq q\leq \infty$. The
modifications in case $q<1$ are straight-forward. Similar as done in the
previous proof we obtain by triangle inequality 
\begin{equation}\label{tri}
   \Big[\sum\limits_{j\in \N_{0}^d}
2^{r|j|_1q}\omega_{m}^{e(j)}(f,2^{-j})^q_p\Big]^{1/q} \lesssim 
   \sum\limits_{\ell \in \zz^d} \Big[\sum\limits_{j\in \N_{0}^d}
2^{r|j|_1q}\omega_{m}^{e(j)}(f_{j+\ell},2^{-j})^q_p\Big]^{1/q}\,,
\end{equation}
where we put (in contrast to the previous proof)
$$
    f_{j}(x) = \sum\limits_{k \in \D_{j}} d^2_{j,k}(f)v_{j,k}(x)\,,
$$
with $f_j = 0$ if $j\notin \N_{-1}^d$\,.
We exploit the piecewise linearity of the basis functions $v_{j+\ell,k}$
together with the at least second order differences in 
$\omega_{m}^{e(j)}(f_{j+\ell},2^{-j})_p$. In fact, let us consider the variable
$x_1$. For $\ell_1<0$ and $|h_1|<2^{-j_1}$ the difference
$\Delta^m_{h_1}(v_{j_1+\ell_1,k_1},x_1)$ vanishes unless $x_1$ belongs to one
of the intervals $I_L, I_M, I_R$ given by $I_L:=\{x\in
\tor:|x-2^{-(j_1+\ell_1)}k_1|\lesssim |h_1|\}$, $I_M:=\{x\in
\tor:|x-2^{-(j_1+\ell_1)}(k_1+1/2)|\lesssim |h_1|\}$, and $I_R:=\{x\in
\tor:|x-2^{-(\ell_1+j_1)}(k_1+1)|\lesssim |h_1|\}$. Furthermore, in case
$\ell_1<0$ it is easy 
to verify that 
$$
    |\Delta^m_{h_1}(v_{j_1+\ell_1,k_1},x_1)| \lesssim 2^{\ell_1} \quad,\quad
x_1\in I_L\cup I_M\cup I_R\,.
$$
 In particular, as a consequence of
$|I_L\cup I_M \cup I_R| \lesssim |h_1| \leq 2^{-j_1}$ we obtain
\begin{equation}\label{f22}
   \int_{\tor} |(\Delta^2_{h_1}v_{j_1+\ell_1,k_1})(x_1)|^p\,\dint x_1 \lesssim
2^{p\ell_1}2^{-j_1}
\end{equation}
in case $\ell_1 < 0$\,. In case $\ell_1\geq 0$ we use the trivial estimate 
\begin{equation}\label{f22b}
   \int_{\tor} |(\Delta^m_{h_1}v_{j_1+\ell_1,k_1})(x_1)|^p\,\dint x_1 \lesssim
2^{-(j_1+\ell_1)}\,.
\end{equation}
Now we combine the component-wise estimates in \eqref{f22} and \eqref{f22b} to
estimate $\omega_{m}^{e(j)}(f_{j+\ell},2^{-j})_p$ from above. 
Indeed, using the perfect localization property of the basis functions we
obtain	
$$
    \omega_{m}^{e(j)}(f_{j+\ell},2^{-j})_p \lesssim \Big(2^{-(j+\ell)}\sum_{k\in
\D_{j+\ell}}|d^2_{j+\ell,k}(f)|^p\Big)^{1/p}\prod\limits_{i=1}^d
    \left\{\begin{array}{rcl}
      2^{\ell_i(1+1/p)}&:&\ell_i < 0\,,\\
      1&:&\ell_i\geq 0\,.
    \end{array}\right.
$$
Now, similar as in \eqref{anorm2} in the previous proof we see for $\ell \in
\zz^d$
\begin{equation}
  \begin{split}
   &\sum\limits_{j\in \N_{0}^d}2^{r|j|_1r
q}\omega_{m}^{e(j)}(f_{j+\ell},2^{-j})^q_p\\
   &~~~\lesssim \sum\limits_{j\in \N_0^d} 2^{(r-1/p)|j+\ell|_1 q}\Big(\sum_{k\in
\D_{j+\ell}}|d^2_{j+\ell,k}(f)|^p\Big)^{q/p}
       \left\{\begin{array}{rcl}
      2^{\ell_i(1+1/p-r)q}&:&\ell_i < 0\,,\\
      2^{-r\ell_i q}&:&\ell_i\geq 0\,,
    \end{array}\right.
  \end{split}  
\end{equation}
which results in 
$$
  \sum\limits_{j\in \N_{0}^d}2^{r|j|_1r
q}\omega_{m}^{e(j)}(f_{j+\ell},2^{-j})^q_p
  \lesssim 2^{-q \delta|\ell|_1}\|d^2_{j,k}(f)|s^r_{p,q}b\|^q\,,
$$
where we used that $0<r<1+1/p$\,.
Plugging this into \eqref{tri} concludes the proof\,. \eproof

\begin{rem}\label{rectmeans} The restriction $p \geq 1$ in Proposition
\ref{contvsdisc} can be removed. Note, that this restriction is caused by the
difference characterization in Lemma \ref{diff} which can be extended to
$0<p,q\leq \infty$ and $m>r>(1/p-1)_+$, see \cite[Theorem 4.5.1]{Ul06}, 
by using rectangle means of differences, 
\begin{equation}\label{rectm}
\mathcal{R}_{m}^e(f,t)_p:=
\Big\|\int_{[-1,1]^d}|\Delta_{(h_1t_1,\ldots,h_dt_d)}^{m,e}(f,\cdot)|\dint
h\Big\|_p\quad,\quad t\in [0,1]^d\,,
\end{equation}
instead of the mixed moduli of continuity \eqref{modc}. In other words, if
$$
    0<p,q\leq \infty\quad\mbox{and}\quad (1/p-1)_+<r<1+1/p
$$
it holds with $m\geq \max\{2,1+1/p\}$
\be\label{f101}
    \|f|S^r_{p,q}B(\tor^d)\|^{\varphi}\lesssim
\|f|S^r_{p,q}B(\tor^d)\|_{\mathcal{R}}^{(m)} \lesssim
\|d^2_{j,k}(f)|s^r_{p,q}b\|
\ee
for all $f\in C(\tor^d)$ with finite discrete quasi-norm
$\|d^2_{j,k}(f)|s^r_{p,q}b\|$\,. The quasi-norm in the middle of \eqref{f101}
represents 
the counterpart of \eqref{f4}, where \eqref{modc} is replaced by
\eqref{rectm}\,. Note, that the restriction $r<2$ does not occur here in case
$p<1$.

\end{rem}

\section{Discrepancy of order $2$ digital nets}
\label{discr}

\subsection{Digital $(t,n,d)$-nets}

We quote from \cite[Section 4]{Di08} to describe the construction method of
order $\sigma$ digital $(t,n,d)$-nets which in case $\sigma = 1$ are original
digital nets from \cite{Ni87} but in this form they were introduced in
\cite{Di07}.

For $s,n\in\N$ with $s\geq n$ let $C_1,\ldots,C_d$ be $s\times n$ matrices over
$\F_2$. For $\nu\in\{0,1,\ldots,2^n-1\}$ with the dyadic expansion $\nu = \nu_0
+ \nu_1 2 + \ldots + \nu_{n-1} 2^{n-1}$ with digits
$\nu_0,\nu_1,\ldots,\nu_{n-1}\in\{0,1\}$ the dyadic digit vector $\bar{\nu}$ is
given as $\bar{\nu} = (\nu_0,\nu_1,\ldots,\nu_{n-1})^{\top}\in\F_2^n$. Then we
compute
$C_i\bar{\nu}=(x_{i,\nu,1},x_{i,\nu,2},\ldots,x_{i,\nu,s})^{\top}\in\F_2^s$ for
$1\leq i\leq d$. Finally we define
\[ x_{i,\nu}=x_{i,\nu,1} 2^{-1}+x_{i,\nu,2} 2^{-2}+\ldots+x_{i,\nu,s} 2^{-s}
\in[0,1) \]
and $x_{\nu}=(x_{1,\nu},\ldots,x_{d,\nu})$. We call the point set
$\P_n=\{x_0,x_1,\ldots,x_{2^n-1}\}$ a digital net over $\F_2$.

Now let $\sigma\in\N$ and suppose $s\geq \sigma n$. Let $0\leq t\leq \sigma n$
be an integer. For every $1\leq i\leq d$ we write
$C_i=(c_{i,1},\ldots,c_{i,s})^{\top}$ where $c_{i,1},\ldots,c_{i,s}\in\F_2^n$
are the row vectors of $C_i$. If for all $1\leq
\lambda_{i,1}<\ldots<\lambda_{i,\eta_i}\leq s,\,1\leq i\leq d$ with
\[
\lambda_{1,1}+\ldots+\lambda_{1,\min\{\eta_1,\sigma\}}+\ldots+\lambda_{d,1}
+\ldots+\lambda_{d,\min\{\eta_d,\sigma\}}\leq\sigma n - t \]
the vectors
$c_{1,\lambda_{1,1}},\ldots,c_{1,\lambda_{1,\eta_1}},\ldots,c_{d,\lambda_{d,1}},
\ldots,c_{d,\lambda_{d,\eta_d}}$ are linearly independent over $\F_2$, then
$\P_n$ is called an order $\sigma$ digital $(t,n,d)$-net over $\F_2$.

The quality parameter $t$ and the order $\sigma$ qualify the structure of the
point set, the lower $t$ and the higher $\sigma$ -- the more structure do the
point sets have.

\begin{lem}
 Let $\P_n$ be an order $1$ digital $(t,n,d)$-net then every dyadic interval of
order $n-t$ contains exactly $2^t$ points of $\P_n$.
\end{lem}

Therefore $(t,n,d)$-nets are also $(t+1,n,d)$-nets and order $\sigma+1$ nets are
also order $\sigma$ nets (with even lower quality parameter). In particular
every point set $\P_n$ constructed with the method described above is at least
an order $\sigma$ digital $(\sigma n,n,d)$-net. We refer to \cite{Di07} and
\cite{Di08} for more details on such relations.

\noindent We need the following fact concerning projections of digital nets
(\cite[Theorem 2]{DiKr09}). 
\begin{lem}\label{Kritzer} Let $\mathcal{P}_n$ be an order $2$ digital
$(t,n,d)$-net. Let further $I_\ell:=\{i_1,\ldots,i_\ell\} \subset \{1,\ldots,d\}$ 
be a fixed set of coordinates. Then the projection $\mathcal{P}_n(I_\ell)
\subset [0,1]^{\ell}$ of the set $\mathcal{P}_n$ on the coordinates 
in $I_\ell$ is an order $2$ digital $(t,n,\ell)$-net.  
\end{lem}

Now we quote explicit constructions of higher order digital nets. We will only
briefly describe the method, for details consult \cite{DiPi14_1},
\cite{DiPi14_2} and \cite{Di14}. The starting point are order $1$ digital
$(t',n,\sigma d)$-nets and the so called digit interlacing composition
\begin{align}
\begin{split}
 \mathscr{D}_\sigma:\, [0,1]^\sigma &\rightarrow [0,1]
 \\
 (x_1,\ldots,x_\sigma)\, &\mapsto \sum_{a=1}^\infty \sum_{r=1}^\sigma \xi_{r,a}
2^{-r-(a-1)\sigma}, \label{eqn_interlacing}
 \end{split}
\end{align}
where $\xi_{r,1},\xi_{r,2},\ldots$ are the digits of the dyadic decomposition of
$x_r$. The digit interlacing is applied component wise on vectors, namely
\begin{align*}
 (x_1,\ldots,x_{\sigma d})\, &\mapsto
(\mathscr{D}_\sigma(x_1,\ldots,x_\sigma),\ldots,\mathscr{D}_\sigma(x_{
(d-1)\sigma+1},\ldots,x_{\sigma d})).
\end{align*}
element wise. 
Suppose that $\P_n$ is an order $1$ digital $(t',n,\sigma d)$-net. Then
$\mathscr{D}_\sigma^d(\P)$ is an order $\sigma$ digital $(t,n,d)$-net with $
t=\sigma t' + d\sigma(\sigma-1)/2$. Therefore, it is possible to construct order
$2$ digital $(t,n,d)$-nets.

\subsection{The discrepancy function and its Haar coefficients}

Let $N$ be a positive integer and let $\P$ be a point set in $[0,1]^d$ with $N$
points. Then the discrepancy function $D_{\P}$ is defined as
\begin{align}
D_{\P}(x) = \frac{1}{N} \sum_{z \in \P} \chi_{[0,x]}(z) - x_1 \cdots x_d
\end{align}
for any $x = (x_1, \ldots, x_d) \in [0,1]^d$. By $\chi_{[0,x]}$ we denote the
characteristic function of the interval $[0,x] =
[0,x_1]\times\ldots\times[0,x_d]$, so the term $\sum_z \chi_{[0,x]}(z)$ is equal
to $\# (\P \cap [0,x])$. $D_{\P}$ measures the deviation of the number of points
of $\P$ in $[0,x]$ from the fair number of points $N |[0,x]| = N \, x_1 \cdots
x_d$.

For further studies of the discrepancy function we refer to the monographs
\cite{DiPi10, NoWo10, Ma99, KuNi74} and surveys \cite{Bi11, Hi14, MaDiss13}.

We will use Haar coefficients of the discrepancy function which are given by
\cite[Proposition 5.7]{Ma15}\,.

\begin{prop}\label{lev14} Let $\P_{n}$ be an order $2$ digital $(t,n,d)$-net
over $\F_2$. Let further $j\in \N_0^d$ and $m\in \D_j$. Then there exists a
constant $c = c(d) > 0$ that satisfies the following properties.\\
{\em (i)} If $|j|_1 \geq n-\lceil t/2 \rceil$ then
$$
      \mu_{j,m}(D_{\mathcal{P}_n})\leq c\, 2^{-|j|_1-n+t/2}
$$
and $\mu_{j,m}(D_{\mathcal{P}_n}) \leq c\, 2^{-2|j|_1}$ for all but $2^n$ values
of $m$\,.\\
{\em (ii)} If $|j|_1 < n-\lceil t/2 \rceil$ then 
$$
      \mu_{j,m}(D_{\mathcal{P}_n}) \leq c \, 2^{-2n+t}(2n-t-2|j|_1)^{d-1}\,.
$$
\end{prop}

\begin{rem}\label{Hinrichs} {\em (i)} It is in this proposition that the higher
order property of the digital nets is needed. For a usual order 1  digital
$(t,n,d)$-net
the main factor in the second estimate would only be $2^{-|j|_1-n-t}$ instead of
$2^{-2n+t}$ which is not sufficient to yield the
right order of convergence in our results.\\
{\em (ii)} In \cite{Hi10} the author computed Haar coefficients of the
dicrepancy function with respect to the two-dimensional Hammersley point set.
In contrast to the current paper the author considered the non-periodic
situation $S^r_{p,q}B([0,1]^d)$ and therefore had to deal with the cases
where components $j_i$ of $j$ equal $-1$ as well. This corresponds to the inner product
of $D_{\mathcal{P}_n}$ with the characteristic function $\chi_{[0,1]}$ in the
respective component $i$, see \cite[Theorem 3.1, (iv)]{Hi10}. This upper bound
is essentially sharp, as shown in \cite[Lemma 3.7]{Hi10}, and responsible
for the fact, that the results in \cite{Hi10} can not be extended to $r<0$. A
$d$-dimensional counterpart for Chen-Skriganov points can be found in
\cite{MaDiss13}.
\end{rem}

\section{QMC integration for periodic mixed Besov spaces}
\label{QMC}

In the sequel we consider quasi-Monte Carlo integration methods for
approximating the integral $I(f):=\int_{\tor^d} f(x)\,\dint x$ of a $d$-variate
 continuous function $f\in C(\tor^d)$. More precisely, for a discrete set
$\mathcal{P} \subset [0,1]^d$
of $N$ points we compute 
$$
  I_N(\mathcal{P},f) := \frac{1}{N}\sum\limits_{z\in \mathcal{P}}
f(z)\quad,\quad f\in F_d\,,
$$
where $F_d$ denotes a class of functions from $C(\tor^d)$. Assume that for $f\in
F_d$ the multivariate Faber expansion 
\eqref{repr} converges in $C(\tor^d)$. We consider the integration error
$R_N(f):=I_N(\mathcal{P},f)-I(f)$\,. In fact,
\begin{equation}\label{interror}
  \begin{split}
   |R_N(f)|  &= \Big|\frac{1}{N}\sum\limits_{z\in \mathcal{P}}f(z) -
\int_{\tor^d}f(x)\,\dint x\Big| \\
  & = \Big|\sum\limits_{j\in \N_{-1}^d} \sum\limits_{m\in \D_j}
d^2_{j,m}(f)\frac{1}{N}\sum\limits_{z\in \mathcal{P}}v_{j,m}(z)-
\sum\limits_{j\in \N_{-1}^d} \sum\limits_{m\in \D_j}
d^2_{j,m}(f)\int_{\tor^d} v_{j,m}(x)\,\dint x\Big|\\
&= \Big|\sum\limits_{j\in \N_{-1}^d} \sum\limits_{m\in \D_j}
d^2_{j,m}(f)c_{j,m}(\mathcal{P})\Big|\,,
  \end{split}
\end{equation}
where
\begin{equation}\label{f32}
   c_{j,m}(\mathcal{P}):=\frac{1}{N}\sum\limits_{z\in \mathcal{P}}v_{j,m}(z)-
   \int_{\tor^d} v_{j,m}(x)\,\dint x\quad,\quad j\in \N_{-1}^d, m\in \D_j\,.
\end{equation}

Let us first take a look at the second summand. 

\begin{lem}\label{intvj} Let $j\in \N_{-1}^d$ and $m\in \D_j$ then
$$
    \int_{\tor^d} v_{j,m}(x)\,\dint x = 2^{-|j+1|_1}\,.
$$
\end{lem}

\bproof We use the tensor product structure of the $v_{j,m}$ to compute
$$
    \int_{\tor^d} v_{j,m}(x)\,\dint x = \prod\limits_{i=1}^d
\int_{\tor}v_{j_i,m_i}(x_i)\,\dint x_i = \prod\limits_{i=1}^d 2^{-|j_i+1|}
    = 2^{-|j+1|_1}\,.
$$
\eproof 

The next Lemma connects the Haar coefficients $\mu_{j,m}(D_{\mathcal{P}})$ of
the discrepancy function $D_{\mathcal{P}}$ with the numbers
$c_{j,m}(\mathcal{P})$\,.

\begin{lem}\label{haarfaber} Let $\mathcal{P} \subset [0,1]^d$ with $\# \P =
N$.\\
{\em (i)} If $j\in \N_0^d$ and $m\in \D_j$ we have
$$
    \mu_{j,m}(D_{\mathcal{P}}) =
(-1)^d2^{-|j|_1}c_{j,m}(\mathcal{P})\quad,\quad 
$$
{\em (ii)} If $j\in \N_{-1}^d\setminus \N_0^d$ and $m\in \D_j$ we have
$$
    \mu_{j,m}(D_{\bar{\mathcal{P}}}) = (-1)^s 2^{-|j|_+}c_{j,m}(\mathcal{P})\,,
$$
where $\bar{\mathcal{P}}$ denotes the projection of $\mathcal{P}$ onto those $s$
coordinates $z_i$ where $j_i \neq -1$. Moreover
$\mu_{j,m}(D_{\bar{\mathcal{P}}})$ is the Haar coefficient with respect to the
$s$-variate function $D_{\bar{\mathcal{P}}}$
  
\end{lem}

\bproof Let $j\in \N_0^d$ and $m\in \D_j$. We compute
$\mu_{j,m}(D_{\mathcal{P}})$. This involves two parts. We first deal with
\begin{equation}\label{haar1}
  \begin{split}
     \int_{[0,1]^d}\frac{1}{N}\sum\limits_{z \in \mathcal{P}}
\chi_{[0,x)}(z)h_{j,m}(x)\,\dint x &= 
     \frac{1}{N}\sum\limits_{z \in \mathcal{P}}\int_{[0,1]^d}
\chi_{[z,1)}(x)h_{j,m}(x)\,\dint x\\
     & = \frac{1}{N}\sum\limits_{z \in \mathcal{P}} \prod\limits_{i=1}^d
\int_{z_i}^1 h_{j_i,m_i}(y)\,\dint y\,.
  \end{split}
\end{equation}
Let us deal with the univariate integrals on the right-hand side of
\eqref{haar1}. Clearly, for any $i=1,\ldots,d$,
\begin{equation}\label{inthaar}
  \begin{split}
     \int_{z_i}^1 h_{j_i,m_i}(y)\,\dint y &= -\int_{0}^{z_i}
h_{j_i,m_i}(y)\,\dint y
     = -\int_{0}^{z_i} h(2^{j_i}y-m_i)\,\dint y\\
     &=-2^{-j_i}\int_0^{2^{j_i}z_i-m_i} h(\tau)\,\dint \tau \\
     &= -2^{-j_i} v(2^{j_i}z_i-m_i) = -2^{-j_1}v_{j_i,m_i}(z_i)\,.
  \end{split}
\end{equation}
This together with \eqref{haar1} yields
\be
     \int_{[0,1]^d}\frac{1}{N}\sum\limits_{z \in \mathcal{P}}
\chi_{[0,x)}(z)h_{j,m}(x)\,\dint x
     = (-1)^d 2^{-|j|_1}\frac{1}{N}\sum\limits_{z \in
\mathcal{P}}v_{j,m}(z)\quad,\quad z\in [0,1]^d\,.
\ee
It remains to compute 
\be\label{inthaar2}
    \int_{[0,1]^d}x_1\cdots x_d \,h_{j,m}(x)\,\dint x =
\prod\limits_{i=1}^d \int_0^1 y\,h_{j_i,m_i}(y)\,\dint y\,.
\ee
Integration by parts together with \eqref{inthaar} yields for $i=1,\ldots,d$
$$
  \int_0^1 y\,h_{j_i,m_i}(y)\,\dint y = -2^{-j_i}\int_0^1 v_{j_i,m_i}(y)\,\dint
y
$$
which, together with \eqref{inthaar2}, implies
\be\label{inthaar3}
   \int_{[0,1]^d}x_1\cdots x_d \,h_{j,m}(x)\,\dint x =
(-1)^d2^{-|j|_1}\int_{[0,1]^d} v_{j,m}(x)\,\dint x
\ee
Combining, \eqref{haar1}, \eqref{inthaar}, and \eqref{inthaar3} yields the
result in (i). The result in (ii) is a simple
consequence of (i).
\eproof

The following result represents our main theorem. 
\begin{satz} \label{main_theorem} Let $\mathcal{P}_n$ be an order $2$ digital
$(t,n,d)$-net over $\F_2$. Then for $1\leq p,q\leq \infty$ and 
$1/p<r<2$ there exists a constant $c = c(p,q,r,d) > 0$ and we have with $N=2^n$
$$
    \Int_N(S^r_{p,q}B(\tor^d))\le \sup\limits_{\|f|S^r_{p,q}B(\tor^d)\|\leq 1}
|I(f) - I_N(\mathcal{P}_n,f)| \leq c\, 2^{rt/2} \, N^{-r}(\log
N)^{(d-1)(1-1/q)}\quad,\quad n\in \N\,.
$$
\end{satz}

\bproof Let $f\in S^r_{p,q}B(\tor^d)$. By the embedding result in Lemma
\ref{emb}/(ii),(iii) we see that $f\in S^{\varepsilon}_{\infty,1}B(\tor^d)$ for
an
$\varepsilon>0$\,. As a consequence of Proposition \ref{discvscont} we obtain
that \eqref{repr} converges to $f$ in $C(\tor^d)$ and therefore in
$L_p(\tor^d)$. Then, by \eqref{interror} together with twice H\"older's
inequality we obtain
\begin{equation}\label{f33}
  \begin{split}
    |R_N(f)| \leq& \sum\limits_{j\in \N_{-1}^d} \sum\limits_{m\in \D_j}
		    |d^2_{j,m}(f)c_{j,m}(\mathcal{P}_n)|\\
    \leq& \Big[\sum\limits_{j\in
\N_{-1}^d}2^{(r-1/p)|j|_1q}\Big(\sum\limits_{m\in
\D_j}|d^2_{j,m}(f)|^p\Big)^{q/p}\Big]^{1/q}\\
&\times \Big[\sum\limits_{j\in \N_{-1}^d}2^{-(r-1/p)|j|_1q'}
\Big(\sum\limits_{m\in
\D_j}|c_{j,m}(\mathcal{P}_n)|^{p'}\Big)^{q'/p'}\Big]^{1/q'}\\
\lesssim& \|f|S^r_{p,q}B(\tor^d)\|\cdot
\Big[\sum\limits_{j\in \N_{-1}^d}2^{-(r-1/p)|j|_1q'}
\Big(\sum\limits_{m\in
\D_j}|c_{j,m}(\mathcal{P}_n)|^{p'}\Big)^{q'/p'}\Big]^{1/q'}\,,
  \end{split}
\end{equation}
where we used Proposition \ref{discvscont} in the last step. In order to prove
the error bound 
it remains to estimate the second factor. Let us deal with 
$$
    \sum\limits_{\substack{|j|_1 < n-\lceil t/2 \rceil\\j\in
\N_0^d}}2^{-(r-1/p)|j|_1q'}
\Big(\sum\limits_{m\in \D_j}|c_{j,m}|^{p'}\Big)^{q'/p'}
$$
first. By Lemma \ref{haarfaber} together with 
Proposition \ref{lev14},(ii) we obtain
\begin{equation}
\begin{split}
  &\sum\limits_{\substack{|j|_1<n-\lceil t/2 \rceil\\j\in
\N_0^d}}2^{-(r-1/p)|j|_1q'}
  \Big(\sum\limits_{m\in \D_j}|c_{j,m}|^{p'}\Big)^{q'/p'}\\ 
  &\lesssim \sum\limits_{\substack{|j|_1<n-\lceil t/2 \rceil\\j\in \N_0^d}}
2^{-|j|_1(r-1/p)q'}2^{(-2n+t)q'}2^{|j|_1q'}2^{|j|_1q'/p'}(2n-t-2|j|_1)^{(d-1)q'}
\\
  &\lesssim 2^{(-2n+t)q'}\sum\limits_{\substack{|j|_1<n-\lceil t/2 \rceil\\j\in
\N_0^d}}
  2^{-|j|_1(r-2)q'}(2n-t-2|j|_1)^{(d-1)q'}\\
  &\asymp 2^{(-2n+t)q'}\sum\limits_{\ell=0}^{n-\lceil t/2 \rceil-1}
\ell^{d-1}(2n-t-2\ell)^{(d-1)q'}2^{-\ell(r-2)q'}\,.
\end{split}  
\end{equation}
Putting $M:=n-\lceil t/2 \rceil$ we obtain 
\begin{equation*}
\begin{split}
\sum\limits_{\substack{|j|_1<n-\lceil t/2 \rceil\\j\in
\N_0^d}}2^{-(r-1/p)|j|_1q'}
  \Big(\sum\limits_{m\in \D_j}|c_{j,m}|^{p'}\Big)^{q'/p'}
  \asymp&~ 2^{(-2n+t)q'}\sum\limits_{\ell=0}^{M-1}
\ell^{d-1}(M-\ell)^{(d-1)q'}2^{-\ell(r-2)q'}\\
  \asymp&~ 2^{(-2n+t)q'}2^{(2-r)Mq'}M^{d-1}\\
  &\times \sum\limits_{\ell=0}^{M-1}
2^{(2-r)(\ell-M)q'}(\ell/M)^{d-1}(M-\ell)^{(d-1)q'}\,.
\end{split}
\end{equation*}
At this point we need the assumption $r<2$ in order to estimate
$$
   \sum\limits_{\ell=0}^{M-1}
2^{(2-r)(\ell-M)q'}(\ell/M)^{d-1}(M-\ell)^{(d-1)q'} \leq
\sum\limits_{\ell=0}^\infty
   2^{-\ell(2-r)q'}\ell^{(d-1)q'}\leq C<\infty\,.
$$
This gives 
$$
    \sum\limits_{\substack{|j|_1<n-\lceil t/2 \rceil\\j\in
\N_0^d}}2^{-(r-1/p)|j|_1q'}
  \Big(\sum\limits_{m\in \D_j}|c_{j,m}|^{p'}\Big)^{q'/p'} \lesssim
2^{(-2n+t)q'}2^{(2-r)Mq'}M^{d-1}\asymp 
  2^{-rnq'}\,2^{rq't/2}\,n^{(d-1)}\,.
$$
Let us now deal with $\sum_{|j|_1\geq n-\lceil t/2\rceil}$\,. By Proposition
4.1, (ii) and Lemma \ref{haarfaber} we get 
\begin{equation*}
 \begin{split}
   &\sum\limits_{\substack{|j|_1\geq n-\lceil t/2 \rceil\\j\in
\N_0^d}}2^{-(r-1/p)|j|_1q'}
  \Big(\sum\limits_{m\in \D_j}|c_{j,m}|^{p'}\Big)^{q'/p'}\\
  &\lesssim \sum\limits_{\substack{|j|_1\geq n-\lceil t/2 \rceil\\j\in
\N_0^d}}2^{-(r-1/p)|j|_1q'}
  \Big(\sum\limits_{m\in A_j}|c_{j,m}|^{p'}\Big)^{q'/p'}
  +\sum\limits_{\substack{|j|_1\geq n-\lceil t/2 \rceil\\j\in
\N_0^d}}2^{-(r-1/p)|j|_1q'}
  \Big(\sum\limits_{m\in \D_j\setminus A_j}|c_{j,m}|^{p'}\Big)^{q'/p'}\,,
 \end{split}
\end{equation*}
where $A_j$ denotes the set of indices $m\in \D_j$ where $I_{j,m} \cap
\mathcal{P}_n \neq \emptyset$. Clearly $\# A_j \leq 2^n$. By Lemma 
\ref{intvj} we directly obtain $|c_{j,m}|\leq 2^{-|j+1|_1}$ if $m\in
D_j\setminus A_j$, whereas by Lemma \ref{haarfaber} and Proposition \ref{lev14},
(i), $|c_{j,m}| \lesssim 2^{-n+t/2}$ if $m\in A_j$. This gives
\begin{equation}
 \begin{split}
      &\sum\limits_{\substack{|j|_1\geq n-\lceil t/2 \rceil\\j\in
\N_0^d}}2^{-(r-1/p)|j|_1q'}
  \Big(\sum\limits_{m\in A_j}|c_{j,m}|^{p'}\Big)^{q'/p'}
  \lesssim 
  \sum\limits_{\substack{|j|_1\geq n-\lceil t/2 \rceil\\j\in \N_0^d}}
2^{-(r-1/p)|j|_1q'}2^{nq'/p'}2^{(-n+t/2)q'}\\
  &~~~~~~~~~~~\lesssim 2^{(-n+t/2)q'}2^{nq'/p'}\sum\limits_{\substack{|j|_1\geq
n-\lceil t/2 \rceil\\j\in \N_0^d}} 2^{-(r-1/p)|j|_1q'}\\
  &~~~~~~~~~~~\lesssim 2^{(-n+t/2)q'}2^{nq'/p'}2^{-(r-1/p)nq'}n^{d-1}\asymp
2^{-rnq'}\, 2^{rq't/2}\,n^{d-1}\,,
 \end{split}
\end{equation} 
where we used $r>1/p$ in the last step.
Furthermore,
\begin{equation*}
 \begin{split}
      &\sum\limits_{\substack{|j|_1\geq n-\lceil t/2 \rceil\\j\in
\N_0^d}}2^{-(r-1/p)|j|_1q'}
  \Big(\sum\limits_{m\in \D_j\setminus A_j}|c_{j,m}|^{p'}\Big)^{q'/p'}
  \lesssim 
  \sum\limits_{\substack{|j|_1\geq n-\lceil t/2 \rceil\\j\in \N_0^d}}
2^{-(r-1/p)|j|_1q'}2^{-|j|_1 q'}2^{|j_1|q'/p'}\\
  &~~~~~~~~~~~~~~~~~\lesssim \sum\limits_{\substack{|j|_1\geq n-\lceil t/2
\rceil\\j\in \N_0^d}}
  2^{-|j|_1(r-1/p+1-1/p')q'} \lesssim 2^{-nrq'}n^{d-1}\,.
 \end{split}
\end{equation*}
Putting everything together yields 
\be\label{final}
   \Big[\sum\limits_{j\in \N_0^d}2^{-(r-1/p)|j|_1q'}
\Big(\sum\limits_{m\in
\D_j}|c_{j,m}(\mathcal{P}_n)|^{p'}\Big)^{q'/p'}\Big]^{1/q'} \lesssim 2^{-rn}\,
2^{rt/2}\,n^{(d-1)(1-1/q)}\,.  
\ee
It remains to consider the sum $\sum_{j\in \N_{-1}^d\setminus \N_0^d}$\,. In
fact, we decompose
$$
    \sum\limits_{j\in \N_{-1}^d} = \sum\limits_{e \subset \{1,\ldots,d\}}
\sum\limits_{j\in \N_{-1}^d(e)}\,,
$$
where $\N_{-1}^d(e) = \{j\in \N_{-1}^d:j_i \in \N_0 \mbox{ if } i \in e, \mbox{
and } j_i = -1 \mbox{ if } i \notin e\}$. By Lemma \ref{haarfaber}, (ii)
together with Lemma \ref{Kritzer} we can estimate $\sum_{j\in \N_{-1}^d(e)}$ by
means of \eqref{final} and obtain for any fixed $e \neq \emptyset$
\be\label{f110}
\Big[\sum\limits_{j\in \N_0^d(e)}2^{-(r-1/p)|j|_1q'}
\Big(\sum\limits_{m\in
\D_j}|c_{j,m}(\mathcal{P}_n)|^{p'}\Big)^{q'/p'}\Big]^{1/q'} \lesssim
2^{-rn}n^{(|e|-1)(1-1/q)}\,.
\ee
Note, that in case $e = \emptyset$ we obtain $c_{(-1,\ldots,-1),0}(\mathcal{P}_n) =
0$\,. Finally, \eqref{final} and \eqref{f110} together yield
$$
   \Big[\sum\limits_{j\in \N_{-1}^d}2^{-(r-1/p)|j|_1q'}
\Big(\sum\limits_{m\in
\D_j}|c_{j,m}(\mathcal{P}_n)|^{p'}\Big)^{q'/p'}\Big]^{1/q'} \lesssim 2^{-rn}\,
2^{rt/2}\,n^{(d-1)(1-1/q)}\,,
$$
which concludes the proof.\eproof

\begin{rem}\label{non-periodic} {\em (i)} We comment on the question of
extending the result in Theorem \ref{main_theorem} to non-periodic Besov-spaces
of mixed smoothness $S^r_{p,q}B([0,1]^d)$. Similar as in \cite{Hi10, MaDiss13}
we are able to prove a corresponding bound for non-periodic spaces if $r\leq
1$ using order $2$ digital nets. By Remark \ref{Hinrichs},(ii) and the
correspondence in Lemma \ref{haarfaber} this result may not extend to $r>1$ in
general. Note, that the Hammersley points represent an order $2$ digital net in
$d=2$. Therefore, the correct order of $\Int_N(S^r_{p,q}B([0,1]^d))$ in the
non-periodic situation $r>1$ is open.\\
{\em (ii)} In order to integrate non-periodic functions one would rather use a
standard modification of the QMC-rule (which then is not a longer a QMC-rule)
given by 
\begin{equation}\label{modcub}
      Q_N^{\psi}(f)\,:=\, I_N(X_N,T^\psi_d f)
			\,=\,\sum\limits_{j=1}^N
\frac{\prod\limits_{i=1}^d \psi'(x^j_i)}{N}\,f(\psi(x^j_1),\ldots,\psi(x^j_d))\,,
\end{equation}
where 
\[
T^{\psi}_d\colon f\mapsto \Big(\prod\limits_{i=1}^d \psi'(x_i)\Big)\cdot
f(\psi(x_1),\ldots,\psi(x_d))\quad,\quad f\in L_1([0,1]^d)\,,x=(x_1,\ldots,x_d)\in
[0,1]^d\,,
\]
and $\psi'$ is a univariate sufficiently smooth bump with $\psi'(0) = \psi'(1)
= 0$. If $T^{\psi}_d:S^r_{p,q}B([0,1]^d) \to S^r_{p,q}B(\T)$ is a bounded
mapping, see \cite{Du2} and the upcoming paper \cite{UU15_2}, then the
cubature formula $Q_N^{\psi}$ has the same order of convergence on
$S^r_{p,q}B([0,1]^d)$ as $I_N(X_N,\cdot)$ on the periodic space
$S^r_{p,q}B(\T)$.

\end{rem}

The result in Theorem \ref{main_theorem} is optimal. In fact, the following
lower bound for general cubature rules has been shown in \cite{Te90} and
independently with a different method in \cite{DiUl14}.

\begin{satz} \label{Theorem[lower bound]}
Let $1 \le p,q \le \infty$ and $r > 1/p$. Then we have
\begin{equation} \nonumber
\mbox{Int}_N(S^r_{p,q}B(\tor^d))
\ \gtrsim \
N^{-r} (\log N)^{(d-1)(1 - 1/q)}\quad,\quad N\in \N\,.
\end{equation}
\end{satz}

By embedding, see Lemma \ref{emb}, (iv), (v) we directly obtain the following
bound for the classes $S^r_pH(\tor^d)$.

\begin{cor}\label{mainH} Let $\mathcal{P}_n$ be an order $2$ digital
$(t,n,d)$-net over $\F_2$. Then for $1<p<\infty$ and $\max\{1/p,1/2\}<r<2$ there
exists a constant $c=c(p,q,r,d)>0$ and we have with $N=2^n$
$$
   \Int_N(S^r_pH(\tor^d)) \leq \sup\limits_{\|f|S^r_pH(\tor^d)\|\leq 1} |I(f) -
I_N(\mathcal{P}_n,f)| \leq c\, 2^{rt/2} \, N^{-r}(\log N)^{(d-1)/2}\quad,\quad
n\in \N\,.
$$
\end{cor}
\bproof If $p>2$ we use the embedding Lemma \ref{emb}, (v) together with Theorem
\ref{main_theorem}. Note, that we need $r>1/2$ here. If $1<p\leq 2$ we use Lemma
\ref{emb}, (v) together with Theorem \ref{main_theorem}.\eproof

\begin{rem}\label{remH} The case $2<p<\infty$ and $1/p<r\leq 1/2$ is not covered
by Corollary \ref{mainH}. This situation is often referred to as the ``Problem
of small smoothness''. It is not known how digital nets (order 1 should be
enough) behave in this situation. Temlyakov \cite{Te91} was the first who
observed an interesting behavior of the asymptotical error for the Fibonacci
cubature rule in the bivariate situation in spaces $S^r_pH(\tor^2)$. Recently,
in \cite{UlUl14} this behavior has been also established for the Frolov method
in the $d$-variate situation. In fact, for spaces $S^r_pH(\tor^d)^{\square}$
with support in the unit cube $Q_d$ it holds for $1<p<\infty$ and $r>1/p$
$$
  \mbox{Int}_N(S^r_pH(Q_d)^{\square})\lesssim N^{-r}\left\{\begin{array}{lcl}
    (\log N)^{(d-1)(1-r)}&:&p>2 \wedge 1/p<r<1/2,\\
    (\log N)^{(d-1)/2}\sqrt{\log\log N}&:&p>2 \wedge r=1/2,\\
    (\log N)^{(d-1)/2}&:&r>\max\{1/p,1/2\}\,.
  \end{array}\right.
$$
We strongly conjecture the same behavior for $S^r_pH(\tor^d)$ where classical
digital (order $1$) nets give the optimal upper bound.

\end{rem}

\begin{figure}[H]
	\includegraphics[width=0.45\linewidth,
height=0.43\linewidth]{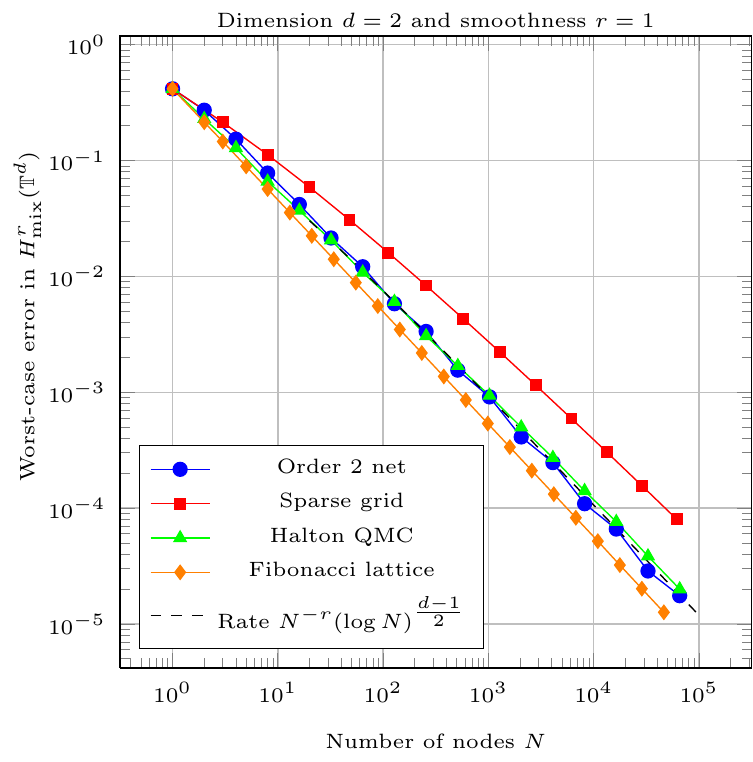}
	\includegraphics[width=0.45\linewidth,
height=0.43\linewidth]{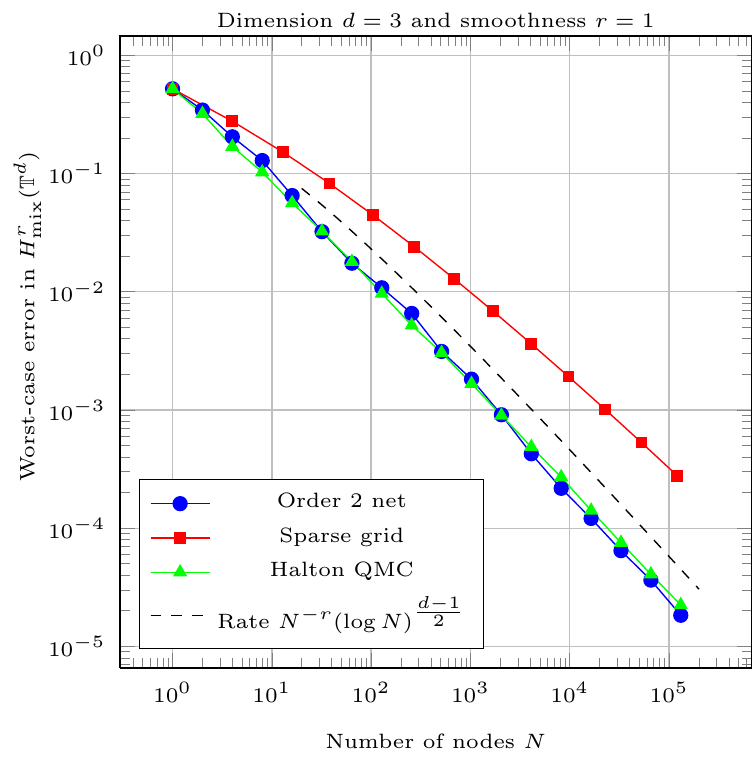}
	\caption{Worst-case errors of order 2 digital nets in
$H^r_\text{mix}(\tor^d)$ for smoothness $r=1$. } \label{fig_errR1_1}
\end{figure}
\section{Numerical experiments}
In this section we use the theory of reproducing kernel Hilbert spaces (RKHS) to
explicitly compute the worst-case error for particular constructions of order 2
digital nets based on Niederreiter-Xing sequences in the case of integer
smoothness $r\in \{1,2\}$. Moreover, we give numerical examples for the case of
fractional smoothness, i.e. $r=\frac{3}{2}$.

\subsection{Worst-case errors in $\mathbf{H^2_{\text{mix}}(\tor^d)}$}
Let us recall that the Besov space $S^r_{2,2} B(\tor^d)$ coincides with the
classical tensor product Sobolev space $H^r_\text{mix}(\tor^d) := H^r(\tor)
\otimes \ldots \otimes H^r(\tor)$ of functions with mixed derivatives of order
$r$ bounded in $L_2(\tor^d)$. Since for $r>1/2$ the space
$H^r_\text{mix}(\tor^d)$ is a Hilbert space which is embedded in $C(\tor^d)$, it
is well-known \cite{Aronszajn_1950} that for a given choice of an inner product
$\langle \cdot, \cdot\rangle_{H^r_{\text{\mix}}}$ there exists a symmetric and
positive definite kernel $K: \tor^d \times \tor^d \rightarrow \re$ that
reproduces point-evaluation, i.e., it holds $f(x) = \langle f(\cdot), K(\cdot,
x) \rangle_{H^r_{\text{\mix}}}$ for all $x \in [0,1)^d$ and $f \in
H^r_{\text{\mix}}(\tor^d)$. Then one can use the well-known worst-case error
formula to compute the quantities
\begin{align}
	\|R_N|(H^{r}_\text{mix})^*\|^2 = & \sup_{\|f|H^r_\text{mix}\| \leq 1}
\left| \int_{\tor^d} f(x) \, d x - \sum_{i=1}^N \lambda_i f(x^i) \right|^2 \\
			    = & \int_{\tor^d} \int_{\tor^d} K(x,y) \, \mathrm{d}
x \mathrm{d} y -2 \sum_{i=1}^N \lambda_i \int_{\tor^d} K(x^i,y) \, \mathrm{d} y
+ \sum_{i=1}^N \sum_{j=1}^{N} \lambda_i \lambda_j K(x^i,x^j)
\end{align}
explicitly, if a point set $\mathcal{P}_N = \{x^1, \ldots, x^N\}$ of 
integration nodes and associated integration weights $\lambda_1, \ldots,
\lambda_N \in \mathbb{R}$ are given. In order to have a simple closed-form
representation of the kernel $K$ we choose the inner product of the univariate
Sobolev space $H^r(\tor)$ to be
\begin{equation}
	\langle f,g\rangle_{H^{r}(\mathbb{T})} = \hat{f}(0) \,
\overline{\hat{g}(0)} + \sum_{k \in \mathbb{Z} \setminus \{0\}} |2\pi k|^{2r}
\hat{f}(k)  \,  \overline{\hat{g}(k)} .
\end{equation}
\begin{figure}[H]
	\includegraphics[width=0.45\linewidth,
height=0.43\linewidth]{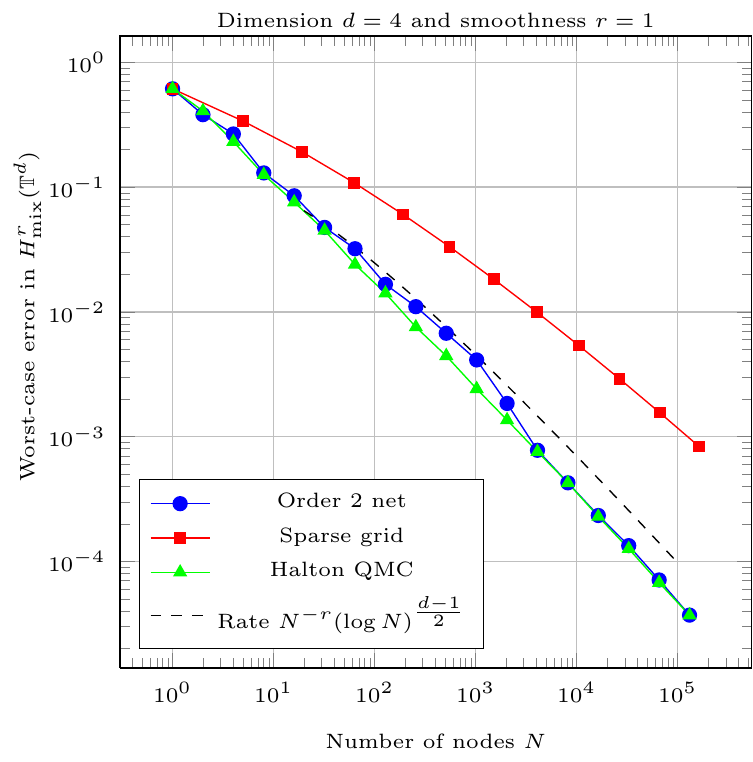}
	\includegraphics[width=0.45\linewidth,
height=0.43\linewidth]{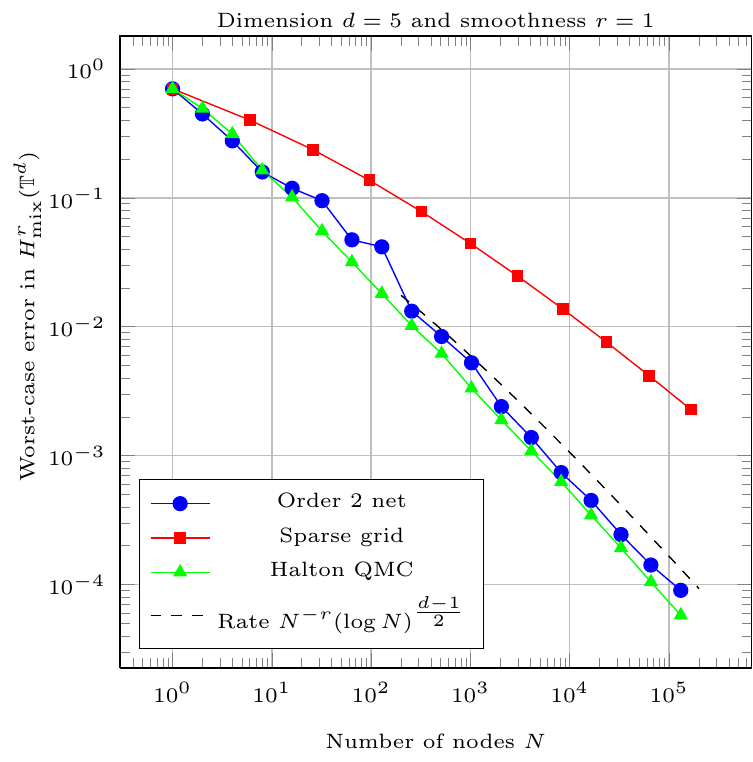}
	\caption{Worst-case errors of order 2 digital nets in
$H^r_\text{mix}(\tor^d)$ for smoothness $r=1$.} \label{fig_errR1_2}
\end{figure}
The induced norm is given by
\begin{equation}
	\|f\|_{H^{r}}^2 = |\hat{f}(0)|^2 + \sum_{k \in \mathbb{Z} \setminus
\{0\}}^\infty |2\pi k|^{2r}  \,  |\hat{f}(k)|^2 =  \left|\int_0^1 f(x) \, d
x\right|^2 + \int_0^1 |f^{(r)}(x)|^2 \, \mathrm{d} x,
\end{equation}
where the last equality only holds for $r \in \mathbb{N}$. Then the reproducing
kernel $K_1: \tor \times \tor \rightarrow \re$ of $H^r(\tor)$ is given by
\cite{Wahba75}
\begin{equation} \label{eqn_repro-kernel}
    \begin{aligned}
	K_{1,r}(x,y) := & 1 + \sum_{k \in \mathbb{Z} \setminus \{0\}} |2\pi
k|^{-2r} \exp(2\pi \mathrm{i} k (x-y)) \\
		   = & 1 + 2\sum_{k=1}^{\infty} |2\pi k|^{-2r} \cos(2\pi k
(x-y))\,.
	\end{aligned}
\end{equation}
If $r \in \mathbb{N}$ the kernel can be written as
\begin{equation} \label{eqn_bernoulli}
	K_{1,r}(x,y) = 1+ \frac{(-1)^{r+1}}{(2r)!} B_{2r}(|x-y|),
\end{equation}
where $B_{2r}: [0,1] \rightarrow \mathbb{R}$ denotes the Bernoulli polynomial of
degree $2r$. Since $H^r_\text{mix}(\tor^d)$ is the tensor product of univariate
Sobolev spaces, the reproducing kernel of $H^r_\text{mix}(\tor^d)$ is given by
the product of the univariate kernels, i.e.
\[
	K_{d,r}(x,y) = \prod_{j=1}^d K_{1,r}(x_j, y_j), \quad x,y \in \tor^d
\]
reproduces point evaluation in $H^r_\text{mix}(\tor^d)$.

As an example we employ order 2 digital nets that are based on Xing-Niederreiter
sequences \cite{NiederreiterXing}, which are known to yield smaller $t$-values
than e.g.\ Sobol- or classical Niederreiter-sequences \cite{Dick2008572}. For
the special case of rational places this construction was implemented by Pirsic
\cite{Pirsic}, 
see also \cite{DiPi10}. It is known \cite{Niederreiter1996241} that one obtains
a digital $(t,n,d)$-net from a digital $(t,n,d-1)$-sequence $\{x^0, \ldots,
x^{2^n-1} \}$ 
by adding an equidistant coordinate, i.e.
\begin{equation} \label{eqn_seq2net}
	\{(\lfloor x^i_1) \rfloor_n, \ldots, \lfloor x^i_{d-1}) \rfloor_n,
i/2^n): i=0,\ldots, 2^n-1  \},
\end{equation}
where $\lfloor\cdot \rfloor_n$ denotes the $n$-th digit truncation.

\begin{figure}[t]
	\includegraphics[width=0.45\linewidth,
height=0.43\linewidth]{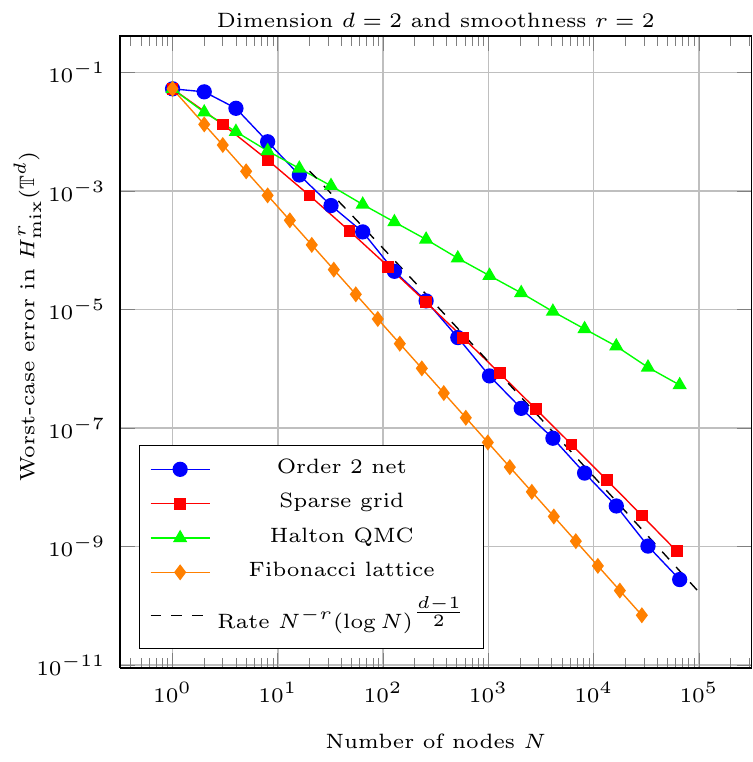}
	\includegraphics[width=0.45\linewidth,
height=0.43\linewidth]{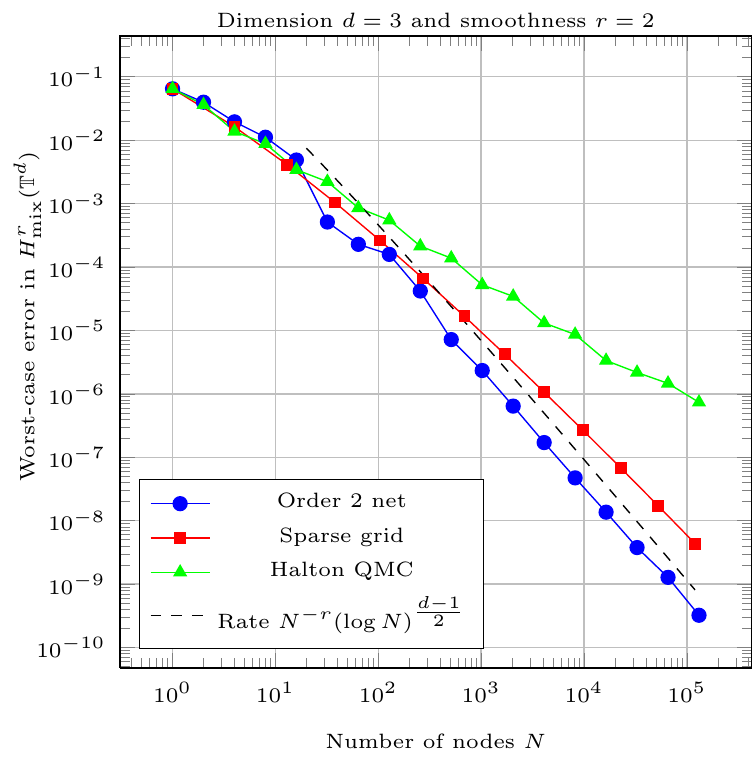}
	\caption{Worst-case errors for smoothness $r=2$ and various cubature
methods.} \label{fig_errR2_1}
\end{figure}

\begin{figure}[t]
	\includegraphics[width=0.45\linewidth,
height=0.43\linewidth]{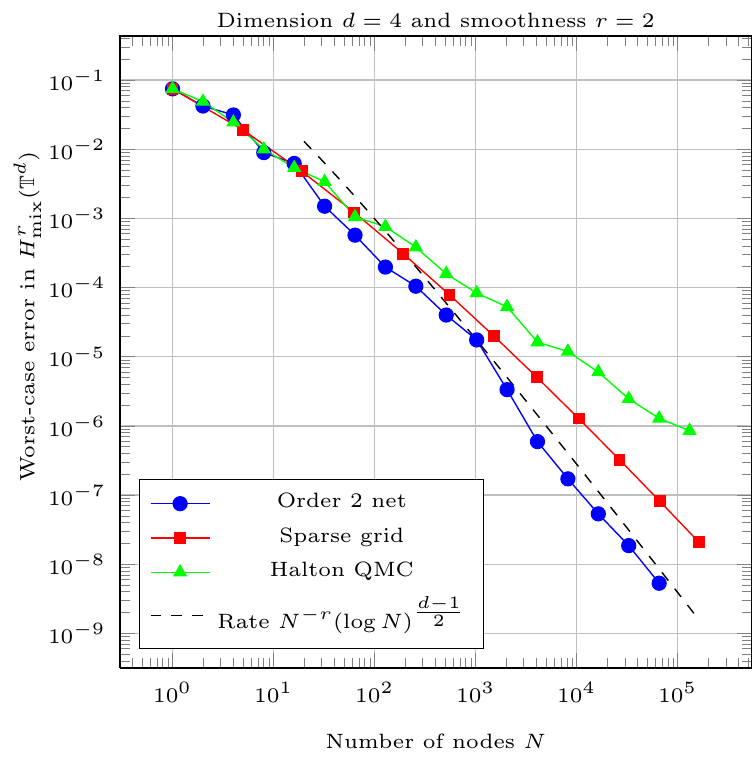}
	\includegraphics[width=0.45\linewidth,
height=0.43\linewidth]{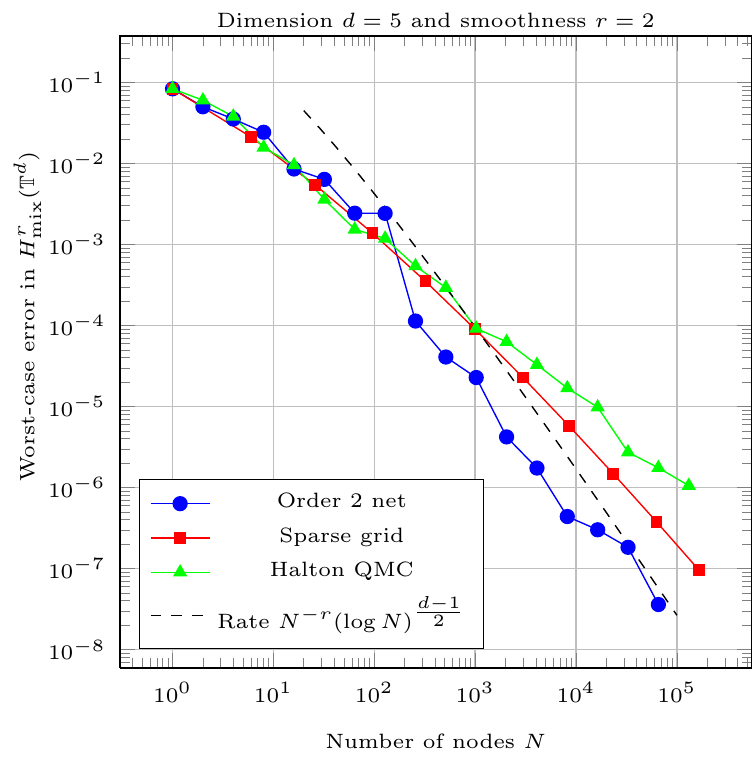}
	\caption{Worst-case errors for smoothness $r=2$ and various cubature
methods.} \label{fig_errR2_2}
\end{figure}

So  we first construct a classical (order 1) digital net from the
Xing-Niederreiter sequence using the 'sequence-to-net' propagation rule
\eqref{eqn_seq2net}. Then we employ the digit interlacing operation
\eqref{eqn_interlacing} to obtain an order 2 net.

For this particular kernel and point set the formula for the squared worst-case
error can be written as
\begin{equation}
	\sup_{\|f|H^r_\text{mix}(\tor^d)\| \leq 1} \left| \int_{\tor^d} f(x) \,
d x - 2^{-n} \sum_{i=0}^{2^n-1}f(x^i) \right|^2 = -1 + 2^{-2n}
\sum_{i=0}^{2^n-1} \sum_{j=0}^{2^n-1} K_{d,r}(x^i, x^j).
\end{equation}
In Figures  \ref{fig_errR1_1} and \ref{fig_errR1_2} we computed the worst-case
errors of the described construction of an order 2 digital net in
$H^1_\text{mix}(\tor^d)$ for dimensions $d=2,\ldots, 5$ and compared it to the
bounds from Theorem \ref{main_theorem}. These expected rates of convergence
$N^{-r}(\log N)^{(d-1)/2}$ were plotted in dashed lines. One can see that our
observed rate of convergence matches the predicted one even though a
dimension-dependent constant seems to be involved. Additionally we computed the
worst-case error for the Halton construction \cite{Halton:1964} which is amongst
the most popular QMC sequences and performs very well for smoothness $r=1$.
Moreover, we consider the sparse grid construction which consists of certain
tensor products of the univariate trapezoidal rule, yielding an error decay of
$\mathcal{O}(N^{-r}(\log N)^{(d-1)(r+\frac{1}{2})})$, see \cite{DiUl14}. Sparse
grids go back to ideas from Smolyak \cite{Smolyak:1963} and belong to today’s
standard approaches when 
it comes 
to high-dimensional problems, see e.g. \cite{Bungartz.Griebel:2004} and the
references therein. It can be seen that their rate of convergence depends
stronger on the dimension $d$ than the low-discrepancy approaches.

The same analysis was done for the case of second order smoothness in 
$H^2_\text{mix}(\tor^d)$. The results are given in Figures \ref{fig_errR2_1} and
\ref{fig_errR2_2}. Here, all quantities were computed in 128-bit floating point
arithmetic.
In the bivariate case it is known \cite{Te91, DiUl14} that the Fibonacci lattice
performs asymptotically optimal. This can also be observed in Figures
\ref{fig_errR1_1} and \ref{fig_errR2_1} where the Fibonacci lattice yields the
same (optimal) rate of convergence as the order 2 digital nets, although it
seems to have a significantly smaller constant. 
For small Fibonacci numbers, it is even known that the Fibonacci lattice is the
globally optimal point set \cite{HO14}.
In summary, we can see that the order $2$ net, the Fibonacci lattice and the
sparse grid are able to benefit from the higher order smoothness, while the
Halton sequence does not improve over $N^{-1}(\log N)^{(d-1)/2}$. 

\subsection{Integration of kink functions}
Mixed Sobolev regularity $H^r_{\text{\mix}}$ is often not suitable to reflect
the correct asymptotical behavior of the integration error of one fixed
function. In case of kink functions, like for instance the Faber hat functions
$v_{j,k}$ from \eqref{f19}, we observe the Sobolev regularity $v_{j,k} \in
H^{3/2-\eps}$ whereas the Besov regularity is $B^2_{1,\infty}$. The tensor
product kink functions belong to $H^{3/2-\eps}_{\text{mix}}$, but as well to
$S^2_{1,\infty}B$. This can be easily deduced from the characterization in Lemma
\ref{diff}. Glancing at Theorem \ref{main_theorem}, we see that the (optimal)
error bound does not depend on the integrability parameter $p$ of the mixed
Besov space $S^r_{p,q}B$. Hence, it seems to be reasonable to ``sacrifice''
integrability in order to gain smoothness which makes our Besov model more
suitable for this issue. Our first example is a typical kink function of the
form $g(x) = \max\{0, h(x)\}$. To be more precise, we consider tensor products
of the univariate (
normalized) function
\begin{equation}\label{kink}
 g(x) = \frac{15 \sqrt{5}}{4} \max\left\{\frac{1}{5} - (x-1/2)^2, 0\right\} ,
\end{equation}
which belongs to $B^2_{1,\infty}(\tor)$ and has integral $\int_0^1 g(x) \, \rd x
= 1$. Hence the tensor 
product function
\begin{equation} 
    g_d(x) := \prod_{j=1}^d g(x_j)\quad,\quad x\in \tor^d\,,
\end{equation}
belongs to $S^2_{1,\infty}B(\tor^d)$ with integral $\int_{\tor^d} g_d(x) \, dx =
1$ and the same holds for the shifted functions
\begin{equation} \label{eqn_testfctB}
    \tilde{g}_d(x, \eta) :=  \prod_{j=1}^d
g(\mathrm{frac}(x_j+\eta_j))\quad,\,\quad x\in \tor^d\, ,
\end{equation}
where $\mathrm{frac}(t) = t - \lfloor t \rfloor$ denotes the fractional
part of $t$ and $\eta \in [0,1]^d$.

In order to obtain smooth convergence rates we compute the maximum error of
$1000$ randomly shifted instances of $\tilde{g}_d$, i.e.
\begin{equation}\label{eqn_testfct}
    \tilde{R}_N(g) = \max_{k \in \{1,\ldots, 1000\}} \left| \int_{\tor^d}
\tilde{g}_d(x, \eta^k) \, \rd x -
    \sum_{i=1}^N \lambda_i\tilde{g}_d(x^i, \eta^k) \right|\, .
\end{equation}

\begin{figure}[t]
    \includegraphics[width=0.45\linewidth,
height=0.43\linewidth]{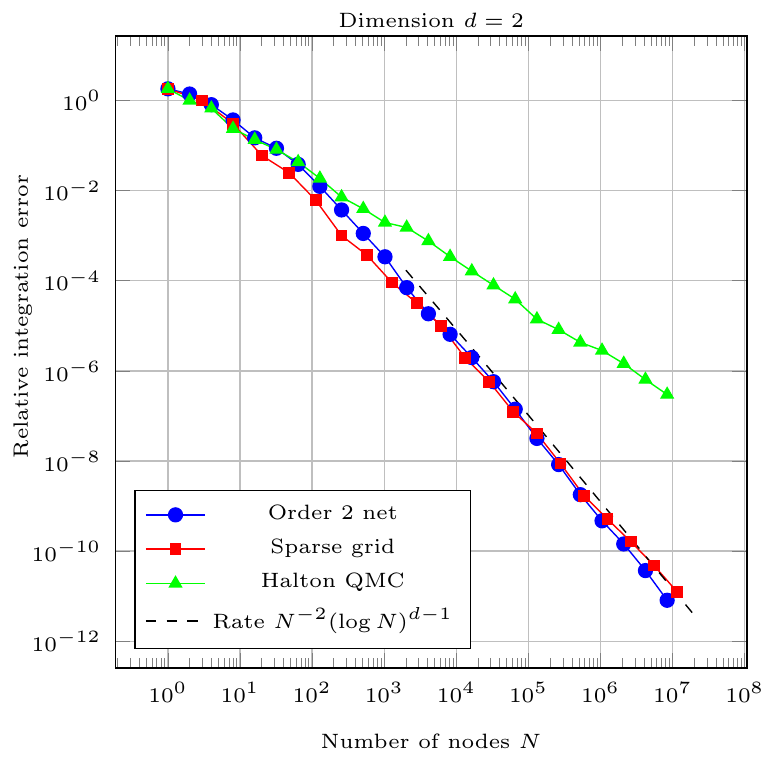}
    \includegraphics[width=0.45\linewidth,
height=0.43\linewidth]{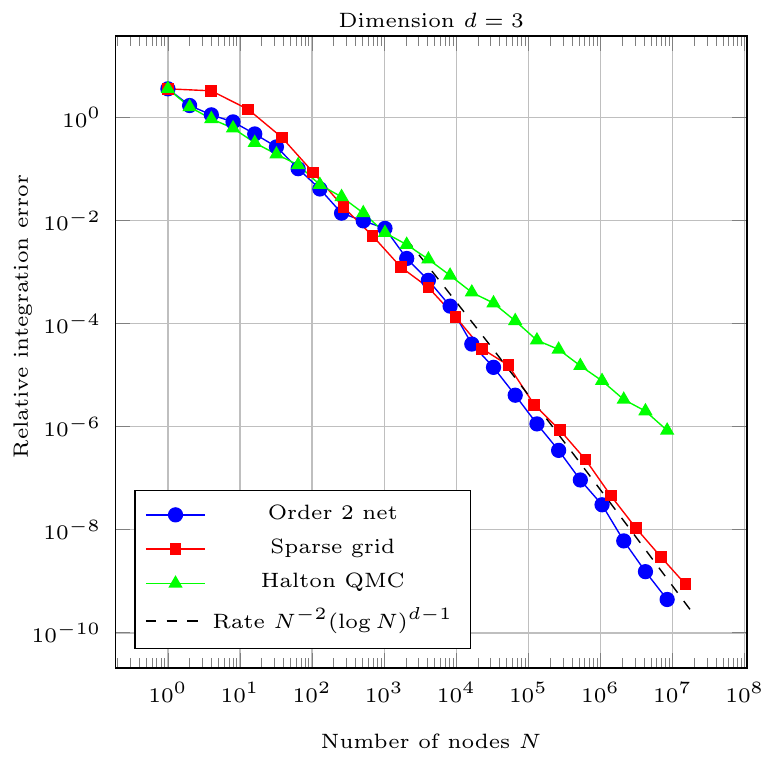}
   
    \includegraphics[width=0.45\linewidth,
height=0.43\linewidth]{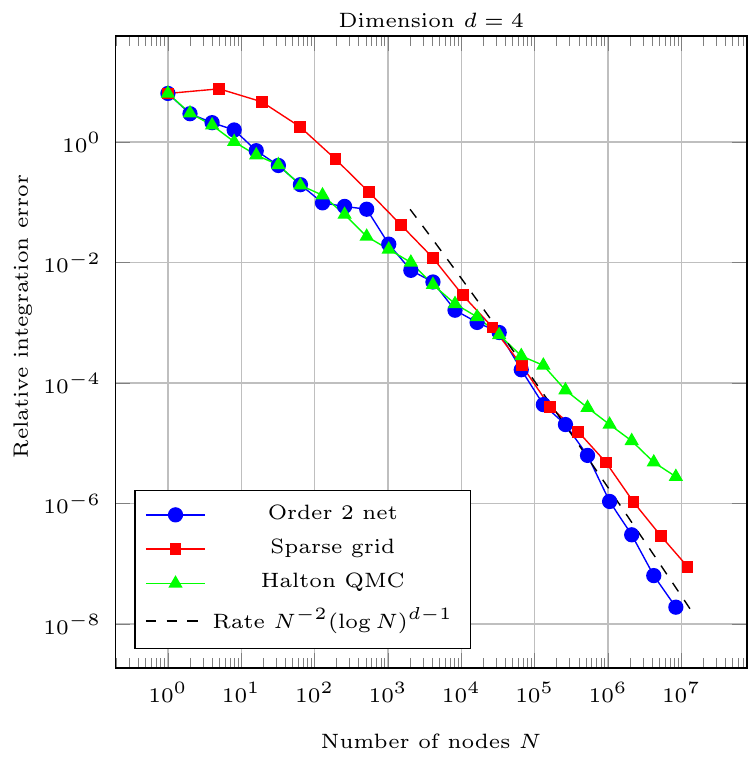}
    \includegraphics[width=0.45\linewidth,
height=0.43\linewidth]{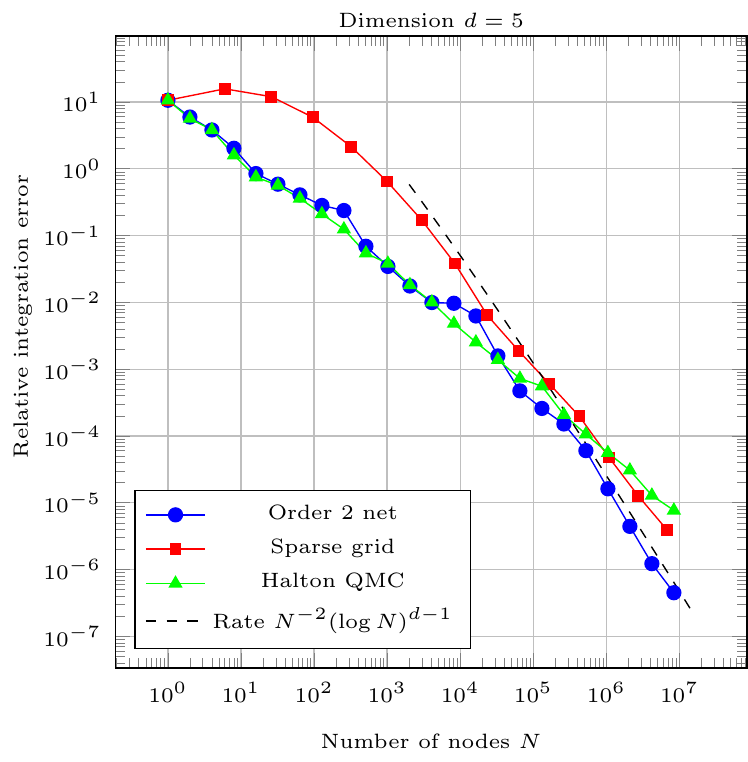}
    \caption{Maximum over the relative integration errors for the test function
\eqref{eqn_testfctB} with kinks.}
\label{fig_errTestFct_B}
\end{figure}

Here, the shifts $\eta^k \sim \mathcal{U}[0,1]^d$
are independent and identically uniformly distributed in $\tor^d$ for
$k=1,\ldots, 1000$ and the integration nodes $x^i, i=1,\ldots,N$ and associated
integration weights $\lambda_i$ depend on the chosen integration method and also
the total number of function values $N$. The results are given in Figure
\ref{fig_errTestFct_B}, where we compared the performance of the order $2$ nets
to both the sparse grid and Halton construction.

\begin{figure}[t]
    \includegraphics[width=0.45\linewidth,
height=0.43\linewidth]{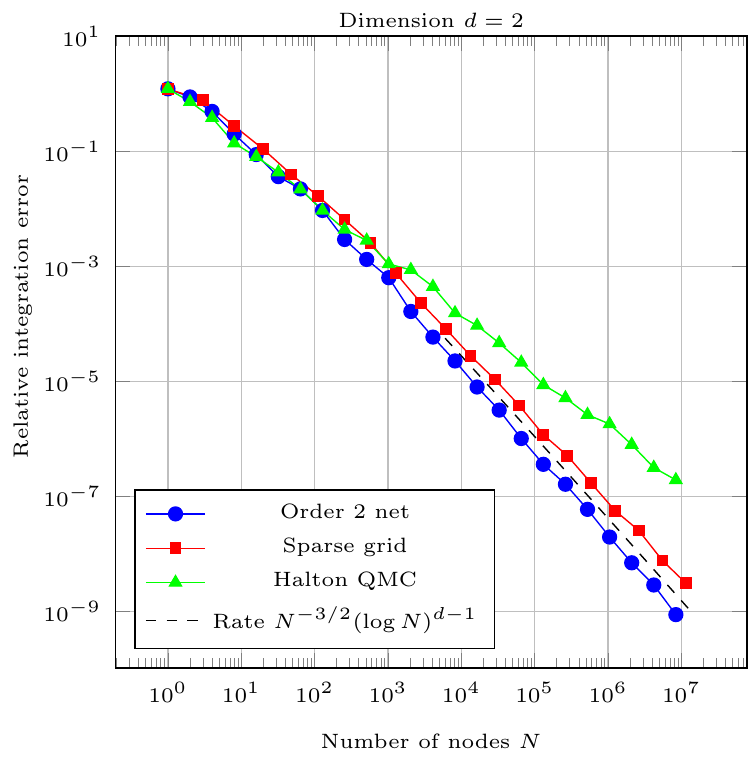}
    \includegraphics[width=0.45\linewidth,
height=0.43\linewidth]{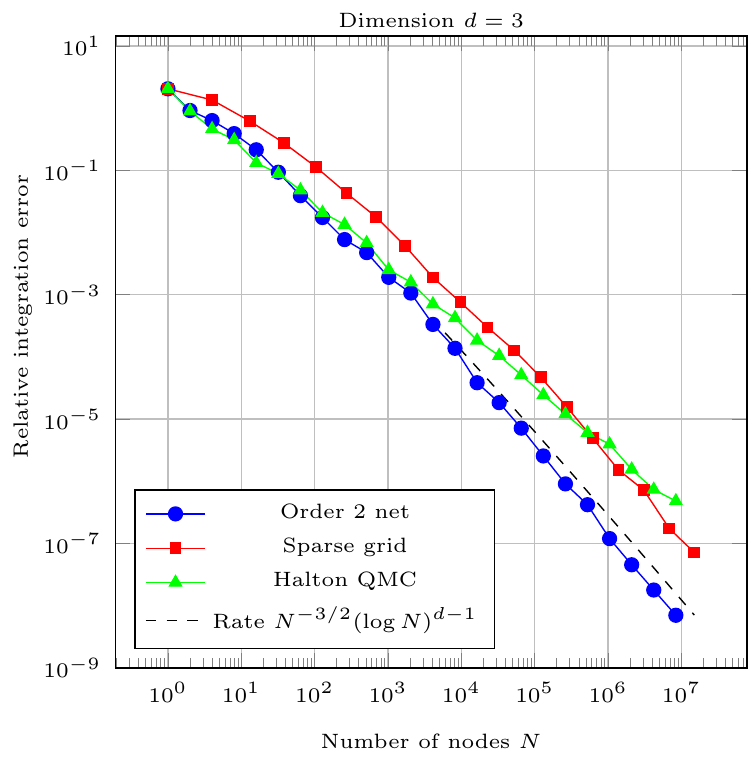}
   
    \includegraphics[width=0.45\linewidth,
height=0.43\linewidth]{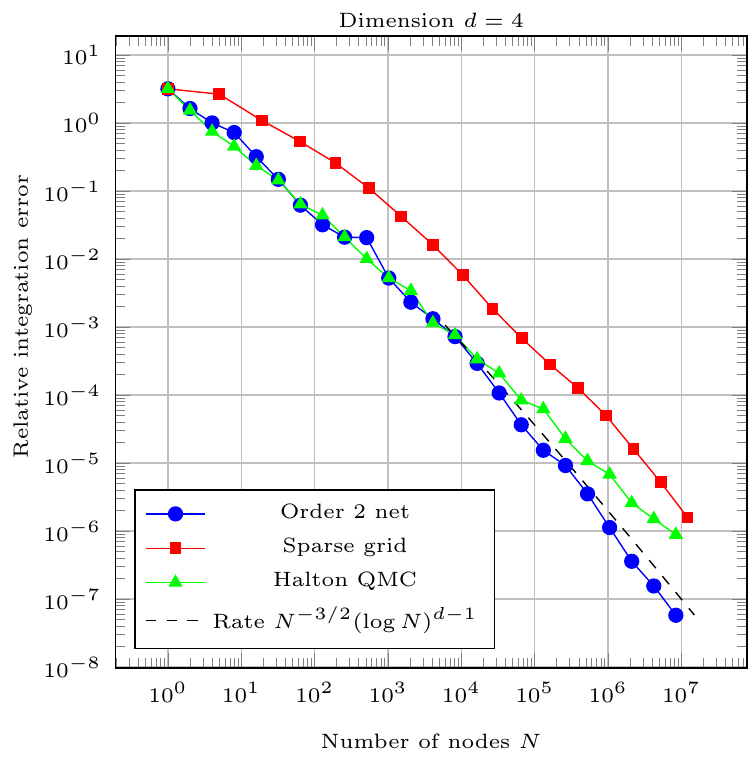}
    \includegraphics[width=0.45\linewidth,
height=0.43\linewidth]{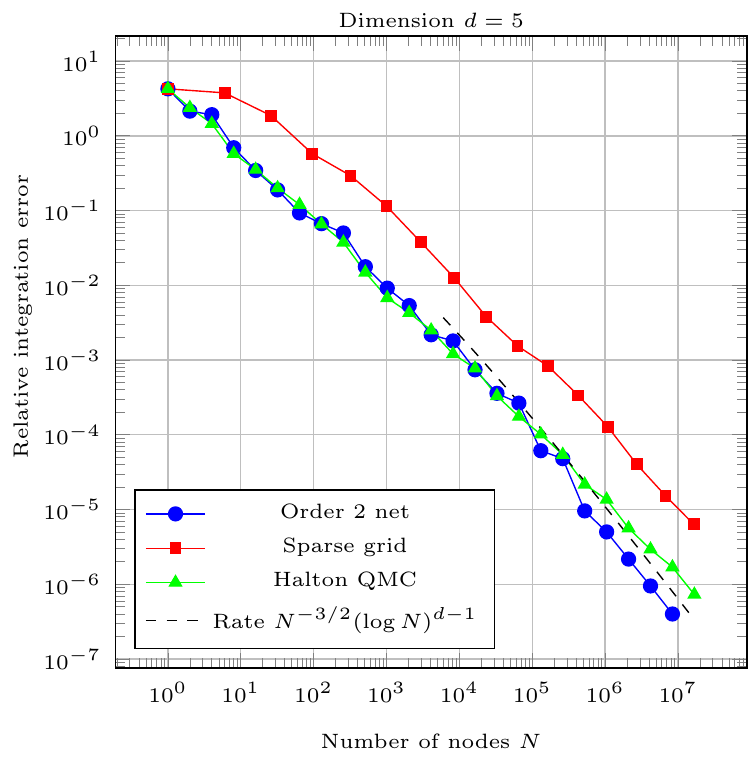}
    \caption{Maximum over the relative integration errors for the test function
\eqref{eqn_testfct_hut} with fractional smoothness $r=\frac{3}{2}$.}
\label{fig_errTestFct_A}
\end{figure}

Next, we consider a toy example from $B^{3/2}_{1,\infty}(\tor)$ which has
Sobolev regularity below $r=1$. We take the square root of the level $0$ hat
function \eqref{f19} normalized with respect to $L_1(\tor^d)$, i.e.,
\begin{equation}    
g(t) := \frac{3}{\sqrt{2}} \sqrt{v_{0,0}(t)}\,.
\end{equation}
It holds $\int_{\tor} g(t) \, d t = 1$. The Besov regularity $r=3/2$ can be
easily deduced from Lemma \ref{diff}. Hence the tensor 
product function 
\begin{equation} 
    g_d(x) := \prod_{j=1}^d g(x_j)\quad,\quad x\in \tor^d\,,
\end{equation}
belongs to $S^{3/2}_{1,\infty}B(\tor^d)$ with integral $\int_{\tor^d} g_d(x) \,
dx = 1$. 
The same holds for the shifted functions
\begin{equation} \label{eqn_testfct_hut}
    \tilde{g}_d(x, \eta) :=  \prod_{j=1}^d
g(\mathrm{frac}(x_j+\eta_j))\quad,\quad x\in \tor^d\,.
\end{equation}
Again, we compute the maximum integration error of $1000$ shifted instances of
$\tilde{g}$. The results are given in Figure \ref{fig_errTestFct_A}. It can be
clearly observed that the obtained convergence rates match the ones predicted in
Theorem \ref{main_theorem}, i.e. $N^{-\frac{3}{2}} (\log N)^{d-1}$.

\noindent

{~}\\
\noindent
{\bf Acknowledgments}
The authors would like to thank the HCM Bonn and the organizers of the
HCM-workshop \emph{Discrepancy, Numerical Integration, and Hyperbolic Cross
Approximation} where this work has been initiated. In addition, they would like
to thank the organizers of the semester program \emph{High-Dimensional
Approximation} at ICERM, Brown University, for providing a pleasant and fruitful
working atmosphere. Finally, they would like to thank Dinh D\~ung, Michael
Griebel, Winfried Sickel and Vladimir Temlyakov for several helpful remarks on
earlier versions of this manuscript.


\begin{thebibliography}{10}

\bibitem{Am76}
T.~I. Amanov.
\newblock {\em Spaces of Differentiable Functions with Dominating Mixed
  Derivatives}.
\newblock Nauka Kaz. SSR, Alma-Ata, 1976.

\bibitem{Aronszajn_1950}
N.~Aronszajn.
\newblock Theory of reproducing kernels.
\newblock {\em Trans. Amer. Math. Soc.}, 68:337--404, 1950.

\bibitem{Ba63}
N.~S. Bakhvalov.
\newblock Optimal convergence bounds for quadrature processes and integration
  methods of {M}onte {C}arlo type for classes of functions.
\newblock {\em Zh. Vychisl. Mat. i Mat. Fiz.}, 4(4):5--63, 1963.

\bibitem{Bi11}
D.~Bilyk.
\newblock On {R}oth's orthogonal function method in discrepancy theory.
\newblock {\em Unif. Distrib. Theory}, (6):143--184, 2011.

\bibitem{Bungartz.Griebel:2004}
H.-J. Bungartz and M.~Griebel.
\newblock Sparse grids.
\newblock {\em Acta Numer.}, 13:1--123, 2004.

\bibitem{Di07}
J.~Dick.
\newblock Explicit constructions of quasi-{M}onte {C}arlo rules for the
  numerical integration of high-dimensional periodic functions.
\newblock {\em SIAM J. Numer. Anal.}, 45:2141--2176, 2007.

\bibitem{Di08}
J.~Dick.
\newblock Walsh spaces containing smooth functions and quasi-{M}onte {C}arlo
  rules of arbitrary high order.
\newblock {\em SIAM J. Numer. Anal.}, 46:1519--1553, 2008.

\bibitem{Di14}
J.~Dick.
\newblock Discrepancy bounds for infinite-dimensional order two digital
  sequences over {$\mathbb{F}_2$}.
\newblock {\em J. Number Theory}, 136:204--232, 2014.

\bibitem{DiKr09}
J.~Dick and P.~Kritzer.
\newblock Duality theory and propagation rules for generalized digital nets.
\newblock {\em Math. Comp.}, 79:993--1017, 2009.

\bibitem{SchwabB}
J.~{Dick}, F.~{Kuo}, Q.~{Thong Le Gia}, and C.~{Schwab}.
\newblock Multi-level higher order qmc galerkin discretization for affine
  parametric operator equations.
\newblock {\em ArXiv e-prints}, June 2014.

\bibitem{SchwabA}
J.~{Dick}, Q.~T. {Le Gia}, and C.~{Schwab}.
\newblock Higher order quasi-monte carlo integration for holomorphic,
  parametric operator equations.
\newblock {\em ArXiv e-prints}, Sept. 2014.

\bibitem{Dick2008572}
J.~Dick and H.~Niederreiter.
\newblock On the exact $t$-value of {N}iederreiter and {S}obol’ sequences.
\newblock {\em J. Complexity}, 24(5–6):572 -- 581, 2008.

\bibitem{DiPi10}
J.~Dick and F.~Pillichshammer.
\newblock {\em Digital nets and sequences. Discrepancy theory and quasi-{M}onte
  {C}arlo integration}.
\newblock Cambridge University Press, Cambridge, 2010.

\bibitem{DP}
J.~Dick and F.~Pillichshammer.
\newblock Discrepancy theory and quasi-monte carlo integration.
\newblock In {\em Panorama in Discrepancy Theory}. Springer--Verlag, 2013.

\bibitem{DiPi14_2}
J.~Dick and F.~Pillichshammer.
\newblock Explicit constructions of point sets and sequences with low
  discrepancy.
\newblock {\em In: Uniform distribution and quasi-Monte Carlo methods -
  Discrepancy, integration and applications}, pages 63--86, 2014.

\bibitem{DiPi14_1}
J.~Dick and F.~Pillichshammer.
\newblock Optimal {$L_2$}-discrepancy bounds for higher order digital sequences
  over the finite field {$\mathbb F_2$}.
\newblock {\em Acta Arith.}, 162(1):65--99, 2014.

\bibitem{Du}
V.~V. Dubinin.
\newblock Cubature formulas for classes of functions with bounded mixed
  difference.
\newblock {\em Matem. USSR Sbornik}, 76:283--292, 1993.

\bibitem{Du2}
V.~V. Dubinin.
\newblock Cubature formulae for {B}esov classes.
\newblock {\em Izvestiya Math.}, 61(2):259--283, 1997.

\bibitem{Di12}
D.~D{\~u}ng.
\newblock B-spline quasi-interpolant representations and sampling recovery of
  functions with mixed smoothness.
\newblock {\em J. Complexity}, 27(6):541--567, 2011.

\bibitem{DiUl14}
D.~D{\~u}ng and T.~{U}llrich.
\newblock Lower bounds for the integration error for multivariate functions
  with mixed smoothness and optimal {F}ibonacci cubature for functions on the
  square.
\newblock {\em Math. Nachr.}, DOI: 10.1002/mana.201400048.

\bibitem{Fa09}
G.~{F}aber.
\newblock {\"U}ber stetige {F}unktionen.
\newblock {\em Math. Ann.}, 66:81--94, 1909.

\bibitem{glasserman2004monte}
P.~Glasserman.
\newblock {\em Monte Carlo Methods in Financial Engineering}.
\newblock Applications of mathematics : stochastic modelling and applied
  probability. Springer, 2004.

\bibitem{Halton:1964}
J.~H. Halton.
\newblock Algorithm 247: Radical-inverse quasi-random point sequence.
\newblock {\em Commun. ACM}, 7(12):701--702, Dec. 1964.

\bibitem{Hi10}
A.~Hinrichs.
\newblock Discrepancy of {H}ammersley points in {B}esov spaces of dominating
  mixed smoothness.
\newblock {\em Math. Nachr.}, 283(3):478--488, 2010.

\bibitem{Hi14}
A.~Hinrichs.
\newblock Discrepancy, integration and tractability.
\newblock {\em In J. Dick, F. Y. Kuo, G. W. Peters, I. H. Sloan, Monte Carlo
  and Quasi-Monte Carlo Methods 2012}, 2014.

\bibitem{HO14}
A.~Hinrichs and J.~Oettershagen.
\newblock Optimal point sets for quasi--{M}onte {C}arlo integration of
  bivariate periodic functions with bounded mixed derivatives.
\newblock to appear in Proc. conf. MCQMC Leuven, 2014.

\bibitem{Hl62}
E.~Hlawka.
\newblock Zur angen\"aherten {B}erechnung mehrfacher {I}ntegrale.
\newblock {\em Monatsh. Math.}, 66:140--151, 1962.

\bibitem{Ko59}
N.~M. Korobov.
\newblock Approximate evaluation of repeated integrals.
\newblock {\em Dokl. Akad. Nauk SSSR}, 124:1207--1210, 1959.

\bibitem{KuNi74}
H.~N. L.~Kuipers.
\newblock {\em Uniform distribution of sequences}.
\newblock John Wiley \& Sons, Ltd., New York, 1974.

\bibitem{MaDiss13}
L.~Markhasin.
\newblock Discrepancy and integration in function spaces with dominating mixed
  smoothness.
\newblock {\em Dissertationes Math.}, 494:1--81, 2013.

\bibitem{Ma13_2}
L.~Markhasin.
\newblock Discrepancy of generalized {H}ammersley type point sets in {B}esov
  spaces with dominating mixed smoothness.
\newblock {\em Unif. Distr. Theory}, 8(1):135--164, 2013.

\bibitem{Ma13_1}
L.~Markhasin.
\newblock Quasi-{M}onte {C}arlo methods for integration of functions with
  dominating mixed smoothness in arbitrary dimension.
\newblock {\em J. Complexity}, 29(5):370--388, 2013.

\bibitem{Ma15}
L.~Markhasin.
\newblock ${L}_2$- and ${S}^r_{p,q}{B}$-discrepancy of (order 2) digital nets.
\newblock {\em Acta Arith.}, 168(2):139--160, 2015.

\bibitem{Ma99}
J.~Matou{\v{s}}ek.
\newblock {\em Geometric discrepancy. An illustrated guide}.
\newblock Springer-Verlag, Berlin, 1999.

\bibitem{UU15_2}
V.~K. Nguyen, M.~Ullrich, and T.~Ullrich.
\newblock Boundedness of pointwise multiplication and change of variable and
  applications to numerical integration.
\newblock {\em Preprint}, 2015.

\bibitem{Ni87}
H.~Niederreiter.
\newblock Point sets and sequences with small discrepancy.
\newblock {\em Monatsh. Math.}, 104:273--337, 1987.

\bibitem{NiederreiterXing}
H.~Niederreiter and C.~Xing.
\newblock A construction of low-discrepancy sequences using global function
  fields.
\newblock {\em Acta Arith.}, 73(1):87--102, 1995.

\bibitem{Niederreiter1996241}
H.~Niederreiter and C.~Xing.
\newblock Low-discrepancy sequences and global function fields with many
  rational places.
\newblock {\em Finite Fields and Their Applications}, 2(3):241 -- 273, 1996.

\bibitem{Nik75}
S.~M. Nikol'skij.
\newblock {\em Approximation of functions of several variables and embedding
  theorems}.
\newblock Nauka Moskva, 1977.

\bibitem{No15}
E.~Novak.
\newblock Some results on the complexity of numerical integration.
\newblock In {\em Monte Carlo and Quasi-Monte Carlo Methods}, to appear in
  Proc. conf. MCQMC Leuven, 2014.

\bibitem{NoWo10}
E.~Novak and H.~Wo{\'z}niakowski.
\newblock {\em Tractability of multivariate problems. {V}olume {II}: {S}tandard
  information for functionals}, volume~12 of {\em EMS Tracts in Mathematics}.
\newblock European Mathematical Society (EMS), Z\"urich, 2010.

\bibitem{Pirsic}
G.~Pirsic.
\newblock {A software implementation of Niederreiter-Xing sequences}.
\newblock In {\em Monte Carlo and quasi-Monte Carlo methods 2000}. Springer,
  Berlin, 2002.

\bibitem{ScTr87}
H.-J. Schmeisser and H.~Triebel.
\newblock {\em Topics in {F}ourier analysis and function spaces}.
\newblock A Wiley-Interscience Publication. John Wiley \& Sons Ltd.,
  Chichester, 1987.

\bibitem{Sk94}
M.~M. Skriganov.
\newblock Constructions of uniform distributions in terms of geometry of
  numbers.
\newblock {\em J. Complexity}, 6:200--230, 1994.

\bibitem{Smolyak:1963}
S.~Smolyak.
\newblock Quadrature and interpolation formulas for tensor products of certain
  classes of functions.
\newblock {\em Dokl. Akad. Nauk SSSR}, 4:240--243, 1963.

\bibitem{Te86}
V.~N. Temlyakov.
\newblock On reconstruction of multivariate periodic functions based on their
  values at the knots of number-theoretical nets.
\newblock {\em Anal. Math.}, 12:287--305, 1986.

\bibitem{Te90}
V.~N. Temlyakov.
\newblock On a way of obtaining lower estimates for the errors of quadrature
  formulas.
\newblock {\em Mat. Sb.}, 181(10):1403--1413, 1990.

\bibitem{Te91}
V.~N. Temlyakov.
\newblock Error estimates for {F}ibonacci quadrature formulas for classes of
  functions with a bounded mixed derivative.
\newblock {\em Trudy Mat. Inst. Steklov.}, 200:327--335, 1991.

\bibitem{Te93}
V.~N. Temlyakov.
\newblock {\em Approximation of periodic functions}.
\newblock Computational Mathematics and Analysis Series. Nova Science
  Publishers Inc., Commack, NY, 1993.

\bibitem{Te03}
V.~N. Temlyakov.
\newblock Cubature formulas, discrepancy, and nonlinear approximation.
\newblock {\em J. Complexity}, 19(3):352--391, 2003.
\newblock Numerical integration and its complexity (Oberwolfach, 2001).

\bibitem{Tr10}
H.~Triebel.
\newblock {\em Bases in function spaces, sampling, discrepancy, numerical
  integration}, volume~11 of {\em EMS Tracts in Mathematics}.
\newblock European Mathematical Society (EMS), Z\"urich, 2010.

\bibitem{Tr12}
H.~Triebel.
\newblock {\em Faber systems and their use in sampling, discrepancy, numerical
  integration}.
\newblock EMS Series of Lectures in Mathematics. European Mathematical Society
  (EMS), Z\"urich, 2012.

\bibitem{UlUl14}
M.~{U}llrich and T.~{U}llrich.
\newblock {T}he role of {F}rolov's cubature formula for functions with bounded
  mixed derivative.
\newblock {\em ArXiv e-prints}, 2015.
\newblock arXiv:1503.08846 [math.NA].

\bibitem{Ul06}
T.~Ullrich.
\newblock Function spaces with dominating mixed smoothness; characterization by
  differences.
\newblock {\em Jenaer Schriften zur Mathematik und Informatik}, Math/Inf/05/06,
  2006.

\bibitem{Ul08}
T.~Ullrich.
\newblock Smolyak's algorithm, sampling on sparse grids and {S}obolev spaces of
  dominating mixed smoothness.
\newblock {\em East J. Approx.}, 14(1):1--38, 2008.

\bibitem{Ul12_3}
T.~{U}llrich.
\newblock Optimal cubature in {B}esov spaces with dominating mixed smoothness
  on the unit square.
\newblock {\em J. Complexity}, 30:72--94, 2014.

\bibitem{Wahba75}
G.~Wahba.
\newblock Smoothing noisy data with spline functions.
\newblock {\em Numer. Math.}, 24(5):383--393, 1975.

\end{thebibliography}
\end{document}